\documentclass[a4paper, 12pt]{article}
\usepackage[utf8]{inputenc}
\usepackage{amsmath,amssymb,amsthm}
\usepackage[a4paper,hmargin=2.5cm,vmargin=2.5cm]{geometry}
\usepackage{tikz-cd} 
\usepackage{xcolor}
\usepackage{graphicx}
\usepackage{cancel}
\usepackage{setspace}
\usepackage{authblk,mathrsfs}
\def\bea{\begin{eqnarray}}
\def\eea{\end{eqnarray}}

\newtheorem{thm}{Theorem}[section]
\newtheorem{cor}[thm]{Corollary}
\newtheorem{ex}[thm]{Example}
\newtheorem{lemma}[thm]{Lemma}

\newtheorem{prop}[thm]{Proposition}
\theoremstyle{definition}
\newtheorem{defn}[thm]{Definition}
\theoremstyle{remark}
\newtheorem{rem}[thm]{Remark}

\DeclareMathOperator{\LHS}{LHS}
\DeclareMathOperator{\RHS}{RHS}

\DeclareMathOperator{\Ker}{Ker}
\DeclareMathOperator{\id}{id}

\DeclareMathOperator{\End}{End}
\DeclareMathOperator{\ad}{ad}
\DeclareMathOperator{\Tr}{Tr}
\DeclareMathOperator{\rk}{rk}
\DeclareMathOperator{\Hom}{Hom}
\DeclareMathOperator{\Mat}{Mat}

\DeclareMathOperator{\Id}{Id}
\DeclareMathOperator{\Fun}{Fun}

	\title{Coboundary ternary structures}
	
	\author[1,2]{M.~A. Chirkov\footnote{mikhlchirkov@gmail.com}}
	\author[1,3,4]{D.~V. Talalaev\footnote{dtalalaev@yandex.ru}}
	\affil[1]{\small Demidov Yaroslavl State University, Yaroslavl, Russia, 150003, Sovetskaya Str. 14
	}
	
	\affil[2]{\small HSE University, Moscow, Russia
	}
	
	\affil[3]{\small Lomonosov Moscow State University, 
		119991, Moscow, Russia.
	}
	
	\affil[4]{\small Kurchatov Institute -- ITEF, Bolshaya Cheremushkinskaya str. 25, 117218, Moscow, Russia.
	}

\begin{document}
\maketitle

\vspace{1em}
\begin{abstract}
Motivated by $r$-matrix Poisson brackets, we study linear, quadratic,
and cubic ternary operations on matrix space determined by an alternating
three-tensor $\varphi$. When these operations satisfy the fundamental
identity, we call them $\varphi$-matrix Nambu brackets. On $\mathfrak{gl}_2$,
the linear operations form a three-dimensional space of Nambu brackets.
For a quadratic Lie algebra with a nonzero central element, we construct
a $3$-Lie bracket from a solution of the modified classical Yang--Baxter
equation and give a rank criterion for the associated linear trivector
to define a Nambu bracket.

We also study ternary operations arising from deformations of (co)associativity. We derive constraints on the deformation coefficients
by expanding the quasi-associative Yang--Baxter equation and the quasi-Hopf
pentagon equation. For a fixed primitive coproduct and a pure exponential
coassociator with an alternating leading tensor, the homogeneous quadratic
pentagon constraint is equivalent to the classical Frenkel--Moore
equation. We show that the Jacobiator of any bivector tangent to a
three-dimensional foliation defines a Nambu bracket. For a cochain twist
on $GL_2$, the bivector defined by the classical limit of the
twisted product has a Nambu Jacobiator whose leaves are the level sets
of the determinant.

\end{abstract}

\clearpage
\tableofcontents

\newpage

\section{Introduction}

A Nambu bracket is a linear map $C^{\infty}(M)^{\wedge n} \to C^{\infty}(M)$ satisfying the fundamental identity (FI, also called the Filippov identity), as well as the Leibniz identity. There is also a similar notion of an $n$-Lie algebra, which differs only in that the $n$-Lie bracket is defined on a vector space and lacks the Leibniz identity. The Nambu bracket is usually denoted in the same way as the Poisson bracket, because when $n = 2$ the fundamental identity coincides with the Jacobi identity. Thus Nambu brackets are a natural generalization of the Poisson bracket to a larger number of functions. It is known that linear Poisson brackets (on a vector space) and Lie algebras are in bijection, but this is not the case for linear Nambu brackets and $n$-Lie algebras: the condition for a Nambu bracket is stricter, namely it requires the local decomposability of the corresponding multivector field where it is nonvanishing.

Originally, a specific Nambu bracket
\begin{equation}\label{canonical-nambu}
    \{f_1, f_2, f_3\} = \det \left( \frac{\partial f_i}{\partial x_j} \right) = \frac{\partial (f_1, f_2, f_3)}{\partial (x_1, x_2, x_3)}
\end{equation}
 on $M = \mathbb{R}^3$ appeared in the work of Y. Nambu \cite{Nambu}. In this work he shows that the vector field $\{\cdot, G, H\}$ defines the Euler equations for the rotation of a rigid body, where $G$ is the total energy of the system and $H$ is the square of the angular momentum, namely
\begin{equation*}
    G = \frac{1}{2} \left( \frac{x^2}{I_x} + \frac{y^2}{I_y} + \frac{z^2}{I_z} \right), H = \frac{1}{2} \left( x^2 + y^2 + z^2 \right)
\end{equation*}
 in coordinates $x = L_1, y = L_2, z = L_3$. Later it became clear that Nambu brackets are applicable in general for the effective description of superintegrable systems, where the flow is given by the vector field $\{\cdot, H_1, \ldots, H_{n-1}\}$, with $H_i$ the integrals of motion \cite{TEMEN2004}. Another useful and rigid property of the Nambu bracket is that fixing $k$ arguments yields a Nambu bracket of arity $n-k$. If $n = 3, k = 1$, then we obtain the Poisson brackets. Moreover, substituting different functions into the Nambu bracket produces compatible Poisson brackets, which are used, for example, for constructing an involutive family of functions.

On the other hand, there is a wide class of Poisson brackets constructed from a solution of the (modified) classical Yang--Baxter equation and a Lie algebra. They are called $r$-matrix Poisson brackets and have developed greatly in connection with Lie bialgebra theory and Poisson-Lie groups \cite{STS}. At the present moment, there is no universally accepted way to generalize these ideas to Nambu brackets. In one paper, a definition of a 3-Lie bialgebra is proposed, but it explicitly uses the presence of an initial 3-Lie algebra \cite{Bai19}. We will mainly work with the associative algebra of matrices, so we do not assume the presence of a natural 3-Lie algebra structure.

Given an alternating tensor $\varphi$, we study the linear operation
$[\varphi,L_1+L_2+L_3]$ and its cubic counterpart
$[\varphi,L_1L_2L_3]$, both of which extend by the Leibniz rule to
triderivations on matrix space.  The central issue is to distinguish
an algebraic 3-Lie identity on linear functions from the
fundamental identity for a Nambu bracket (which is the same for Poisson brackets, but differs for Nambu brackets).
The more general deformation-theoretic question is whether suitable
coherence conditions on a non-associative product imply that
fundamental identity without postulating it separately.

In Section~\ref{sec:nonassoc-rtt}, we study first-order deformations of
(co)associativity and expand the quasi-associative Yang--Baxter equation.
For a fixed primitive coproduct, we derive the constraints imposed by
the pentagon equation at orders $\hbar$ and $\hbar^2$. For an
exponential coassociator with an alternating leading tensor, we show
that the quadratic constraint is equivalent to the classical
Frenkel--Moore equation. Section~\ref{sec:preliminaries} reviews
Nambu brackets, $3$-Lie algebras, and $r$-matrix Poisson brackets.

In Section~\ref{sec:linear-nambu}, we study the general properties of
linear $\varphi$-matrix Nambu brackets and show that on $\mathfrak{gl}_2$
they form a three-dimensional space.
Section~\ref{cubic-section} introduces a cubic analogue of the Sklyanin
bracket and identifies the associated quadratic operation in its expansion
at the identity. In Section~\ref{sec:cyb-fi}, we construct a $3$-Lie
bracket on a reductive Lie algebra equipped with a nondegenerate
invariant symmetric bilinear form. The construction uses a central
element and a solution of the modified classical Yang--Baxter equation.
Also we give a rank criterion for the resulting trivector to define a Nambu bracket.

In Section~\ref{sec:principal-examples}, we compute the Jacobiator of
the bivector arising from a Heisenberg cochain twist. We then show that
the Jacobiator of any bivector tangent to a three-dimensional foliation
defines a Nambu bracket and apply this result to a cochain twist on $GL_2$.

\section{Non-associative algebras and ternary structures}\label{sec:nonassoc-rtt}

A cubic matrix operation related to the construction in Section~\ref{cubic-section} also appears through a deformation
of the RTT relations \cite{ReshetikhinTakhtajanFaddeev} in which associativity is also allowed to fail. This construction is more general because its defining relations do not, by themselves, imply the fundamental identity.
We first describe the resulting compatibility equations.
Let $\mathcal A$ be a (not necessarily associative) unital algebra,
whose generators $t_{ij}$ are assembled into the matrix
\begin{equation*}
 T=\sum_{i,j}E_{ij}\otimes t_{ij}.
\end{equation*}
Consider the matrix relations
\begin{align}
 A T_1(T_2T_3)&=(T_1T_2)T_3 A,\label{quasi-assoc}\\
 R T_1T_2&=T_2T_1R,\label{quasi-comm}
\end{align}
where $A\in\operatorname{End}(V^{\otimes3})$ and
$R\in\operatorname{End}(V^{\otimes2})$ are invertible
operators, with $V=\mathbb C^N$.

In addition to \eqref{quasi-assoc} and \eqref{quasi-comm}, we impose the
quasi-associative Yang--Baxter equation (QAYBE) as an operator identity
on $A$ and $R$ in tensor degree three:
\begin{equation}\label{QAYBE}
 A_{321}^{-1}R_{23}A_{231}R_{13}A_{213}^{-1}R_{12}
 =R_{12}A_{312}^{-1}R_{13}A_{132}R_{23}A_{123}^{-1}.
\end{equation}

For $A=1$, equation \eqref{QAYBE} is the quantum Yang--Baxter
equation.

\subsection{First-order deformation of associativity}
\label{subsubsec:general-deformation}

Consider a flat family $\mathcal A_{\hbar}$ deforming a commutative associative algebra with formal expansions
\begin{equation}\label{AR-first-order}
 R=1+\hbar r+\hbar^2s+O(\hbar^3),\qquad
 A=1+\hbar a+\hbar^2b+O(\hbar^3).
\end{equation}
The first coefficients of the commutator and the associator define
binary and ternary operations. We use the convention
\begin{equation*}
 \operatorname{As}_{\star}(x,y,z)
 =(x\star y)\star z-x\star(y\star z).
\end{equation*}
Assume that, after choosing a vector-space identification of the
formal deformation with $B[[\hbar]]$, where $B$ is a commutative associative algebra,
the product has the form
\begin{equation}\label{star-first}
 x\star y=xy+\hbar B_1(x,y)+\hbar^2B_2(x,y)+O(\hbar^3).
\end{equation}
Then
\begin{align*}
 [x,y]_{\star}&=\hbar\{x,y\}+O(\hbar^2),
 &\{x,y\}&=B_1(x,y)-B_1(y,x),\\
 \operatorname{As}_{\star}(x,y,z)
 &=\hbar\mathfrak a_1(x,y,z)+O(\hbar^2),
\end{align*}
where an elementary calculation gives
\begin{equation*}
 \mathfrak a_1(x,y,z)=B_1(x,y)z+B_1(xy,z)
 -xB_1(y,z)-B_1(x,yz).
\end{equation*}
Thus $\mathfrak a_1$ is a Hochschild coboundary (up to sign). Since $B$ is commutative,
\begin{equation}\label{alt-first-ass}
 \operatorname{Alt}\mathfrak a_1=0.
\end{equation}
If $B_1$ is a
biderivation, then $\mathfrak a_1=0$ without alternation.

Expanding the defining RTT relations \eqref{quasi-assoc}, \eqref{quasi-comm} we have following matrix expressions for binary and ternary operations,
\begin{align}
 \{T_1,T_2\}&=[T_1T_2,r],\label{Skl}\\
 \mathfrak a_1(T_1,T_2,T_3)&=[a,T_1T_2T_3].\label{CubTern}
\end{align}
Since
\eqref{alt-first-ass} holds, a fully alternating first-order defect
must vanish. Thus a nonzero Nambu bracket cannot arise as this
first-order defect of a flat deformation of a commutative algebra.

\subsection{The quasi-associative Yang--Baxter expansion}
\label{subsubsec:qaybe-expansion}

Now we expand \eqref{QAYBE} to see constraints that we have on binary and ternary operation \eqref{Skl}, \eqref{CubTern}. The first-order constraint is independent of $r$ and follows by direct expansion.
\begin{prop}\label{prop:qaybe-first}
If $A$ and $R$ in \eqref{AR-first-order} satisfy QAYBE \eqref{QAYBE}, then
\begin{equation}\label{qaybe-first-obstruction}
 -a_{321}+a_{231}-a_{213}+a_{312}-a_{132}+a_{123}
 =6\operatorname{Alt}(a)=0.
\end{equation}
\end{prop}
This is something really different from Section \ref{cubic-section}, where we have $a \in \mathfrak{g}^{\wedge 3}$.
For the second-order equation, list the first-order coefficients of the six ordered factors on both sides of \eqref{QAYBE}:
\begin{equation}\label{qaybe-sextuples}
\begin{aligned}
 (X_1,\ldots,X_6)
 &=(-a_{321},\ r_{23},\ a_{231},\ r_{13},\ -a_{213},\ r_{12}),\\
 (Y_1,\ldots,Y_6)
 &=(r_{12},\ -a_{312},\ r_{13},\ a_{132},\ r_{23},\ -a_{123}).
\end{aligned}
\end{equation}
Define the following noncommutative quadratic polynomial:
\begin{equation}\label{qaybe-quadratic}
 \mathcal Q_{\mathrm{QYB}}(r,a)
 =\sum_{i<j}X_iX_j-\sum_{i<j}Y_iY_j
 +a_{321}^2+a_{213}^2-a_{312}^2-a_{123}^2.
\end{equation}
\begin{prop}\label{prop:qaybe-second}
If $A$ and $R$ in \eqref{AR-first-order} satisfy \eqref{QAYBE}, then
\begin{equation}\label{qaybe-second}
 6\operatorname{Alt}(b)+\mathcal Q_{\mathrm{QYB}}(r,a)=0.
\end{equation}
\end{prop}
\begin{proof}
Each non-inverse factor has expansion $1+\hbar x+\hbar^2y$,
whereas an inverse factor expands as
$(1+\hbar a+\hbar^2b)^{-1}
=1-\hbar a+\hbar^2(a^2-b)+O(\hbar^3)$.
For six ordered factors with first coefficients $X_i$, the mixed
second-order terms are $\sum_{i<j}X_iX_j$. The four inverse factors
on the two sides contribute the four explicit squares in
\eqref{qaybe-quadratic}. The remaining linear contributions of
$b$ give $6\operatorname{Alt}(b)$, while those of $s$ cancel.
\end{proof}

\begin{cor}\label{cor:qaybe-skew-second}
For $A=1+\hbar^2\varphi+O(\hbar^3)$,
$R=1+\hbar r+O(\hbar^2)$, and
$\varphi\in\wedge^3\mathfrak{gl}_N$, the order-$\hbar^2$
part of \eqref{QAYBE} is
\begin{equation}\label{qaybe-cyb-defect}
 \operatorname{CYB}(r):=[r_{12},r_{13}]+[r_{12},r_{23}]
 +[r_{13},r_{23}]=6\varphi.
\end{equation}
\end{cor}
\begin{proof}
If $a=0$, then \eqref{qaybe-quadratic} is
$-\operatorname{CYB}(r)$. Moreover,
$6\operatorname{Alt}(b)=6\varphi$.
\end{proof}

\subsection{The pentagon equation}
\label{subsubsec:pentagon-expansion}

In this Subsection we consider associative bialgebra, but we deform coproduct in a quasi-coassociative way (analogue of \eqref{quasi-assoc} in coproduct terms). The quasi-bialgebra Pentagon equation compares two paths of reassociation of four tensor factors. To write this equation, suppose that the operator $A$ comes from a universal element, denoted by the same letter, in $H^{\otimes3}[[\hbar]]$.
Here $H$ is an associative unital algebra and
$\Delta:H\to H\otimes H$ is an algebra homomorphism.
Put $A=\Phi^{-1}$, where $\Phi$ satisfies quasi-coassociative condition
\begin{equation}\label{quasi-coassoc-general}
 (\mathrm{id}\otimes\Delta)\Delta(h)
 =\Phi(\Delta\otimes\mathrm{id})\Delta(h)\Phi^{-1}.
\end{equation}
In this convention the pentagon equation is \cite{DrinfeldQuasiHopf,MylonasSchuppSzabo}
\begin{equation}\label{pentagon-A}
\begin{aligned}
 &(\Delta\otimes\mathrm{id}\otimes\mathrm{id})(A)
 (\mathrm{id}\otimes\mathrm{id}\otimes\Delta)(A)\\
 &\qquad=(A\otimes1)
 (\mathrm{id}\otimes\Delta\otimes\mathrm{id})(A)
 (1\otimes A),
\end{aligned}
\end{equation}
which we also abbreviate as
\begin{equation}\label{pentagon-legs}
 A_{12,3,4}A_{1,2,34}=A_{1,2,3}A_{1,23,4}A_{2,3,4}.
\end{equation}
Grouped indices denote coproduct insertions. For example,
$A_{12,3,4}=(\Delta\otimes\id\otimes\id)(A)$.

For the expansion, fix a bialgebra coproduct $\Delta_0$ independent
of $\hbar$ and write
\begin{equation*}
 A=1+\hbar a+\hbar^2b+O(\hbar^3).
\end{equation*}
Introduce the degree-three co-Hochschild coboundary (also called the Cartier coboundary, see \cite{DoninShniderCohom})
\begin{equation}\label{pentagon-D}
 D(u)=u_{12,3,4}+u_{1,2,34}-u_{1,2,3}-u_{1,23,4}-u_{2,3,4},
\end{equation}
where every grouped insertion is defined using $\Delta_0$.
Then the pentagon gives the following exact first two conditions.
\begin{prop}\label{prop:pentagon-two-orders}
For fixed $\Delta = \Delta_0$, equation \eqref{pentagon-A} implies
\begin{align*}
 D(a)&=0,\\
 D(b)+Q(a)&=0,
\end{align*}
where
\begin{equation*}
 Q(a)={}a_{12,3,4}a_{1,2,34}
 -a_{1,2,3}a_{1,23,4}-a_{1,2,3}a_{2,3,4}
 -a_{1,23,4}a_{2,3,4}.
\end{equation*}
\end{prop}
\begin{proof}
Multiply both sides of \eqref{pentagon-legs} formally and compare
coefficients of $\hbar$ and $\hbar^2$. At order two, there is one
cross term on the left and three ordered cross terms on the right.
\end{proof}

For $H=U(\mathfrak g)$ with the primitive coproduct
$\Delta_0(x)=x\otimes1+1\otimes x$, every
$a\in\mathfrak g^{\otimes3}$ satisfies $D(a)=0$.
Co-Hochschild cohomology in this case is
$H^3_{\mathrm{coHoch}}(U(\mathfrak g))\cong
\wedge^3\mathfrak g$ \cite{DoninShniderCohom}. This explains why the alternating
component represents the infinitesimal cohomology class,
although the tensor $a$ need not itself be
alternating.

The exponential ansatz makes the first nonlinear constraint
particularly transparent. Assume that $\Delta=\Delta_0$ is primitive,
$a\in\mathfrak g^{\otimes3}$, and
\begin{equation}\label{pentagon-exponential}
	A=\exp(\hbar a).
\end{equation}
Put $p=a_{123}$, $q=a_{124}$, $u=a_{134}$,
$v=a_{234}$. Since the coproduct is primitive,
\eqref{pentagon-legs} becomes
\begin{equation}\label{pentagon-exponentials}
	e^{\hbar(u+v)}e^{\hbar(p+q)}
	=e^{\hbar p}e^{\hbar(q+u)}e^{\hbar v}.
\end{equation}
\begin{prop}\label{prop:classical-pentagon}
	The coefficient of $\hbar^2$ in \eqref{pentagon-exponentials}
	vanishes if and only if the following homogeneous quadratic
	relation holds:
	\begin{equation}\label{classical-pentagon}
		\begin{aligned}
				\mathcal P(a):={}&[a_{123},a_{124}]
				+2[a_{123},a_{134}]+2[a_{123},a_{234}]\\
				&+[a_{124},a_{134}]+2[a_{124},a_{234}]
				+[a_{134},a_{234}]=0.
		\end{aligned}
	\end{equation}
\end{prop}
\begin{proof}
	The linear terms in \eqref{pentagon-exponentials} agree:
	$(u+v)+(p+q)=p+(q+u)+v$. The order-two difference of
	its left and right sides is one half of
	$[u+v,p+q]-[p,q+u]-[p,v]-[q+u,v]$.
	Expanding the commutators gives $-\mathcal P(a)/2$.
\end{proof}

\subsection{Connection with the Frenkel--Moore equation}
\label{subsubsec:pentagon-FM}

The indices in \eqref{classical-pentagon} suggest a
comparison with the edge-indexed Tetrahedron equation
studied by Frenkel and Moore~\cite{FrenkelMoore1991},
\begin{equation}\label{FM-four-leg}
	\mathscr T_{123}\mathscr T_{124}\mathscr T_{134}\mathscr T_{234}
	=\mathscr T_{234}\mathscr T_{134}\mathscr T_{124}\mathscr T_{123}.
\end{equation}
This version is to be distinguished from Zamolodchikov's tetrahedron
equation with indices $(123),(145),(246),(356)$, they are not equivalent. For
$\mathscr T=1+\hbar a+O(\hbar^2)$, its coefficient of
$\hbar^2$ is the homogeneous quadratic equation
\begin{equation}\label{FM-classical}
	\begin{aligned}
		\mathcal C(a):={}&[a_{123},a_{124}]+[a_{123},a_{134}]
		+[a_{123},a_{234}]\\
		&+[a_{124},a_{134}]+[a_{124},a_{234}]
		+[a_{134},a_{234}]=0.
	\end{aligned}
\end{equation}
Although the coefficients of $\mathcal P$ and $\mathcal C$ are
different, alternation makes the two equations equivalent.

\begin{prop}\label{prop:pentagon-FM-equivalence}
	Let $\mathfrak g$ be a Lie algebra over a field of characteristic
	zero and let $a\in\wedge^3\mathfrak g$, then
	\[
	\mathcal P(a)=0\quad\Longleftrightarrow\quad
		\mathcal C(a)=0.
	\]
\end{prop}
\begin{proof}
	Write $p=a_{123}$, $q=a_{124}$, $u=a_{134}$,
	$v=a_{234}$ and denote by $\sigma = (1243) \in S_4$.
	As $a$ is alternating, $\sigma$ acts by
	\begin{equation}\label{FM-leg-action}
		(p,q,u,v)\longmapsto(q,-v,-p,u).
	\end{equation}
	Direct substitution gives
	\begin{align*}
		\sigma\mathcal P={}&2[p,q]-[p,u]-[p,v]
		+2[q,u]-[q,v]+2[u,v],\\
		\sigma\mathcal C={}&[p,q]-[p,u]-[p,v]
		+[q,u]-[q,v]+[u,v].
	\end{align*}
	Consequently,
	\begin{equation*}
		3\mathcal P+\sigma\mathcal P=5\mathcal C,
			\qquad 2\mathcal P=3\mathcal C-\sigma\mathcal C.
	\end{equation*}
	The vanishing of either expression implies the vanishing
	of its permuted copy. The two displayed identities now prove
	both implications.
\end{proof}

\subsection{Quasi-Hopf bialgebras and ternary structures}\label{subsec:quasi-hopf-background}

The simultaneous appearance of a binary bracket and a ternary operation
in Subsection~\ref{subsubsec:general-deformation} suggests comparison with quasi-Hopf algebras. We
emphasize that three notions must be kept distinct: the coassociator
of a quasi-Hopf coproduct, the associator of a possibly nonassociative star product, and the multiplicative
trivector on a Lie group. They are related in examples, but their
identification is not automatic.

Let $H$ be a quasi-bialgebra with coassociator
$\Phi\in H^{\otimes3}$, characterized, with a fixed convention, by \eqref{quasi-coassoc-general}.
The pentagon identity \eqref{pentagon-A} is an essential additional axiom. Now we describe the setting of the quasi-Hopf deformation in
\cite{MylonasSchuppSzabo}. The
multiplication in $H$ remains associative. In contrast, the
multiplication of an $H$-module algebra may be quasi-associative.
For example, a normalized invertible cochain
$F\in U(\mathfrak g)^{\otimes2}[[\hbar]]$ transforms the
ordinary coproduct into
$\Delta_F(h)=F\Delta(h)F^{-1}$ and yields a coassociator
$\Phi_F$. With a compatible module action the corresponding star
product is
\begin{equation}\label{twist-star}
	f\star g=m\bigl(F^{-1}\triangleright(f\otimes g)\bigr).
\end{equation}
If $F$ fails to be a Hopf 2-cocycle, this star-product need not be
associative, the failure is controlled by the action of $\Phi_F$.

For the expansion \eqref{star-first}, assume that $B_1$ is a
biderivation and put
\begin{equation}\label{binary-classical}
	\{f,g\}:=B_1(f,g)-B_1(g,f)=\pi(df,dg).
\end{equation}
Write $\operatorname{As}_\star=\hbar^2
\operatorname{As}_2+O(\hbar^3)$, because the linear term is zero due to biderivation property.
Then the exact algebraic relation between the Jacobiator of the
commutator and the associator gives
\begin{equation}\label{jacobi-associator}
	J_\pi(f,g,h):=\{f,\{g,h\}\}
		+\{g,\{h,f\}\}+\{h,\{f,g\}\}
		=-6\,\operatorname{Alt}\operatorname{As}_2(f,g,h).
\end{equation}
The right side is therefore a
three-derivation whenever $B_1$ is a biderivation. Note that it need not
satisfy the Fundamental identity.

We denote by $\Pi_J$ the trivector defined by
$\Pi_J(df,dg,dh)=J_\pi(f,g,h)$, which is equal to $\frac{1}{2} [[\pi,\pi]]$. The Schouten Jacobi
identity implies
\begin{equation}\label{schouten-compatibility}
	[[\pi,\Pi_J]] =0.
\end{equation}

\section{Nambu brackets and 3-Lie algebras}\label{sec:preliminaries}
\subsection{Fundamental identity}

Let us define the fundamental identity for a Nambu bracket. For this, it is convenient to introduce the operators $\Delta_{f_1, \ldots, f_{n-1}}: C^{\infty}(M) \to C^{\infty}(M)$ for arbitrary $f_i \in C^{\infty}(M)$.
\begin{align*}
    \Delta_{f_1, \ldots, f_{n-1}}(g) = \{f_1, \ldots, f_{n-1}, g\}.
\end{align*}
Then the fundamental identity can be written as follows:
\begin{equation}\label{fund_id}
    \Delta(\{g_1, \ldots, g_n\}) = \{\Delta(g_1), g_2, \ldots, g_n\} + \ldots + \{g_1, \ldots, g_{n-1}, \Delta(g_n)\},
\end{equation}
where $\Delta = \Delta_{f_1, \ldots, f_{n-1}}$. Even more briefly, the fundamental identity can be written as follows: $\Delta$ is a derivation of the algebra of functions on the manifold, hence it defines a vector field $X = X_{f_1, \ldots, f_{n-1}}$, namely $\mathcal{L}_{X}(g) = \Delta(g)$, where $\mathcal{L}$ is the Lie derivative along the vector field. Such vector fields are called Hamiltonian, by analogy with the Poisson bracket case. Then the fundamental identity means that the flow along any Hamiltonian vector field preserves the Nambu bracket, i.e.
\begin{equation}
    \mathcal{L}_{X} \pi = 0,
\end{equation}
where $\pi$ is the multivector field defining the Nambu bracket. In addition, the Nambu bracket satisfies the Leibniz identity:

\begin{equation}\label{leibniz}
    \{f g, h_1, \ldots, h_{n-1}\} = f \{g, h_1, \ldots, h_{n-1}\} + \{f, h_1, \ldots, h_{n-1}\} g.
\end{equation}

\begin{defn}
    A Nambu bracket (of order $n$) on a smooth manifold $M$ is a multilinear skew-symmetric operation $C^{\infty}(M)^{\wedge n} \to C^{\infty}(M)$ satisfying the fundamental identity \eqref{fund_id} and the Leibniz identity \eqref{leibniz}.
\end{defn}

From now on, we will mainly be interested in the first generalization of the Poisson bracket, namely Nambu brackets for $n = 3$. Let us write the fundamental identity explicitly:
\begin{equation}\label{fund_id_3}
    \{f_1, f_2, \{g_1, g_2, g_3\}\} = \{\{f_1, f_2, g_1\}, g_2, g_3\} + \{g_1, \{f_1, f_2, g_2\}, g_3\} + \{g_1, g_2, \{f_1, f_2, g_3\}\}.
\end{equation}

We shall also speak of a Nambu bracket as a multivector field on $M$, namely $\pi \in \Gamma(\wedge^n T M)$ such that $\{f_1, \ldots, f_n\} = \pi(df_1, \ldots, df_n)$, where $\Gamma(E)$ denotes the space of global sections of the vector bundle $E$. Moreover, if the manifold $M$ is equipped with a volume form $\Omega$, then we are interested in the dual form $\omega \in \Omega_{m-n}(M), \omega = \iota_{\pi} \Omega,$ where $m = \dim M$.

The fundamental identity can be reformulated in the language of multivectors and the dual form as follows.

\begin{thm}[Theorems 2.1 and 2.4 \cite{Vaisman99}]\label{thm-vaisman}
    Let $\pi \in \Gamma(\wedge^n T M)$, $n \geq 3$. Then $\pi$ satisfies the fundamental identity if and only if, at each regular point $(\pi(p) \ne 0)$, the following two conditions hold:
    \begin{itemize}
        \item for any indices $a_1, \ldots, a_{n-1}, b_1, \ldots, b_n$ one has
    \begin{equation*}
        \sum\limits_{u = 1}^{m} \left( \pi^{a_1, \ldots, a_{n-1}, u} \partial_u \pi^{b_1, \ldots, b_n} - \sum\limits_{k = 1}^{n} \pi^{b_1, \ldots, b_{k-1}, u, b_{k+1}, \ldots, b_{n}} \partial_u \pi^{a_1, \ldots, a_{n-1}, b_k} \right) = 0,
    \end{equation*}
        \item $\pi(p)$ is decomposable, equivalently it satisfies the algebraic Pl\"ucker relations.
    \end{itemize}
\end{thm}

The decomposability condition actually arises when the fundamental identity is considered for $f_1 = x_u x_v$ and the remaining functions are coordinate functions. The first condition is exactly the fundamental identity for the functions $f_i = x_{a_i}, g_i = x_{b_i}$. Thus, unlike the Poisson bracket case, where the Jacobi identity is enough to check only on coordinate (linear) functions, for the fundamental identity an additional condition on quadratic functions arises.

\begin{thm}[Theorem 2 \cite{Nakanishi98}]
    Let $\omega = \iota_{\pi} \Omega \in \Omega_{m-n}(M), m > n \geq 3$. Then $\pi$ satisfies the fundamental identity if and only if the following two conditions hold in the neighborhood of each regular point ($\omega(p) \ne 0$):
    \begin{itemize}
        \item there exists a local 1-form $\theta$ such that $d \omega = \theta \wedge \omega$,
        \item there exist local 1-forms $\omega_i$ such that $\omega = \omega_1 \wedge \ldots \wedge \omega_{m-n}$.
    \end{itemize}
\end{thm}

We will also need a useful corollary of this theorem, which is stated without proof in Nakanishi's paper \cite{Nakanishi98}. Since we will use it explicitly, it is good practice to provide a proof.

\begin{cor}[Corollary 2 \cite{Nakanishi98}]\label{corollary_fund}
    Let $\omega = \iota_{\pi} \Omega$, $m = n + 1 = \dim M$. Under these assumptions, $\pi$ satisfies the fundamental identity if and only if
    \begin{equation}\label{fund_condition}
        \omega \wedge d \omega = 0.
    \end{equation}
\end{cor}
\begin{proof}
    The 1-form $\omega$ is automatically decomposable. Near a point where $\omega\ne 0$, complete $\omega$ to a local coframe and write $d\omega=\theta\wedge\omega+\nu$, where $\nu$ involves only the other coframe elements. Wedge multiplication by $\omega$ is injective on these terms, so $\omega\wedge d\omega=0$ if and only if $\nu=0$. The preceding theorem therefore gives the equivalence on the regular locus. On the interior of the zero set both identities hold, since $\omega$ and $\pi$ vanish there. At boundary points the identities extend by smoothness.
\end{proof}

Write $\iota_{df}\pi=\pi(\cdot,\ldots,\cdot,df)$ the multivector field corresponding to fixing the last argument in the Nambu bracket $\{\cdot, \ldots, \cdot, f\}$.

\begin{prop}[Remark 1 \cite{Takhtajan94}]\label{substitution}
    Let $\pi \in \Gamma(\wedge^n T M)$ be a Nambu bracket of order $n$. Then $\iota_{df}\pi \in \Gamma(\wedge^{n-1} T M)$ is a Nambu bracket of order $n-1$.
\end{prop}

\subsection{Relation between Nambu brackets and $n$-Lie algebras}

\begin{defn}[\cite{Filippov86}]
    An $n$-Lie algebra $A$ is a vector space with a completely skew-symmetric multilinear operation $A^{\wedge n} \to A$ satisfying the fundamental identity for any $v_1, \ldots, v_{n-1}, w_1, \ldots, w_n \in A$
\begin{align*}
    [v_1, \ldots, v_{n-1}, [w_1, \ldots, w_n]] &= [[v_1, \ldots, v_{n-1}, w_1], w_2, \ldots, w_n] + \\ &+ \ldots + [w_1, \ldots, w_{n-1}, [v_1, \ldots, v_{n-1}, w_n]].
\end{align*}
\end{defn}

\begin{prop}\label{n-lie-algebra}
    Let $V = \mathbb{R}^m$, $M = V^{*}$, and let $\pi$ be a linear Nambu bracket on $M$, i.e. the coefficients $\pi^{i_1, \ldots, i_n}$ are linear functions of the Cartesian coordinates. Then the Nambu bracket $\pi$ defines an $n$-Lie algebra structure on $V$.
\end{prop}
\begin{proof}
    This is immediate: the coefficients of the linear multivector $\pi$ are precisely structure constants of a skew-symmetric operation on $V$, and the fundamental identity for the Nambu bracket is exactly the fundamental identity for the $n$-Lie algebra. The proof is a direct coordinate check.
\end{proof}

\begin{prop}[Remark 4 in \cite{Takhtajan94}]
    At each regular point of a Nambu bracket of order $n \geq 3$, there is a neighborhood with vector fields $X_1, \ldots, X_n$ such that
    \begin{equation*}
    \pi = X_1 \wedge \ldots \wedge X_n.
    \end{equation*}
\end{prop}

\begin{prop}
    Let $V = \mathbb{R}^m$ and $\pi$ be a linear Nambu bracket on $V^*$. Then the corresponding $n$-Lie algebra has the structure constants determined by the coefficients of $\pi$.
\end{prop}

Now we construct a nontrivial example showing that not every $n$-Lie algebra comes from a decomposable Nambu bracket.

\begin{prop}\label{counterexample}
    For $n = 3$ and $m = 8$, there exists a 3-Lie algebra structure on $\mathbb{R}^8$ whose associated trivector field is not locally decomposable. Hence there are 3-Lie algebras that do not come from a Nambu bracket.
\end{prop}
\begin{proof}
    Consider the 3-Lie algebra $L$ on $\mathbb{R}^4$ with basis $e_1, \ldots, e_4$ and bracket defined by
    \begin{equation*}
        [e_i, e_j, e_k] = \varepsilon^l_{i,j,k} e_l,
    \end{equation*}
    where $\varepsilon^l_{i,j,k} = \varepsilon_{i,j,k,l}$ is the Levi-Civita symbol. Let us denote this 3-Lie algebra by $L$. Then consider the direct sum $L \oplus L$, where the 3-Lie structure is defined independently on each component and all other structural constants vanish. By the standard construction, from the structure constants we obtain a trivector field $\pi = \pi_{(1)} + \pi_{(2)}$ on $\mathbb{R}^8$. We show that it is not decomposable. It is enough to find a point at which the trivector is nonvanishing and not decomposable. Let
    \begin{equation*}
        p = (0,0,0,1,0,0,0,1).
    \end{equation*}
    Then
    \begin{equation*}
        \pi(p) = \partial_1 \wedge \partial_2 \wedge \partial_3 + \partial_5 \wedge \partial_6 \wedge \partial_7.
    \end{equation*}
    How to see that the multivector is not decomposable? Consider the operator $P: \wedge^{n-1} V^* \to V$, $\alpha \mapsto \iota_{\alpha} \pi$. For a decomposable $n$-multivector, we get $\rk P \leq n$. Therefore it is enough to compute $\rk P$. In this case we obtain $\rk P = 6$, hence the multivector $\pi(p)$ is not decomposable.
\end{proof}

\subsection{Lie bialgebras and Poisson-Lie groups}

\begin{defn}[Definition 4.8 \cite{STS}]
    A Lie bialgebra $\mathfrak{g}$ is a vector space with a Lie algebra structure on $\mathfrak{g}$ and on $\mathfrak{g}^{*}$, such that these structures are compatible, namely
    \begin{equation*}
        \delta([X, Y]) = (\ad_X \otimes 1 + 1 \otimes \ad_X)\delta(Y) - (\ad_Y \otimes 1 + 1 \otimes \ad_Y)\delta(X),
    \end{equation*}
    where $\delta: \mathfrak{g} \to \mathfrak{g} \wedge \mathfrak{g}$ is the linear map dual to the Lie bracket on $\mathfrak{g}^*$.
\end{defn}

\begin{defn}[Definition 4.3 \cite{STS}]
    A Poisson-Lie group $G$ is a Lie group with a compatible Poisson bracket, i.e. the multiplication map on the Lie group is a Poisson map $G \times G \to G$.
\end{defn}

\begin{prop}[Proposition 4.7 \cite{STS}]
    The tangent space at the identity of a Poisson-Lie group is a Lie bialgebra.
\end{prop}

How can one construct a compatible Lie bracket on $\mathfrak{g}^*$ if $\mathfrak{g}$ is a Lie algebra? For this, it is useful to use the concept of Lie algebra cohomology.

\begin{defn}[Lie Algebra Cohomology \cite{ChevalleyEilenberg}]
    The cohomology of a Lie algebra $\mathfrak{g}$ with coefficients in a module $M$ is the cohomology of the cochain complex $(C^{\bullet}, d)$, where $C^k(M) = \Hom(\wedge^k \mathfrak{g}, M)$ and $d_k: C^{k} \to C^{k+1}$
    \begin{align*}
        (df)(x_1, \ldots, x_{k+1}) &= \sum_{i} (-1)^{i+1} x_i f(x_1, \ldots, \hat{x}_i, \ldots, x_{k+1}) \, + \\ &+ \sum_{i < j} (-1)^{i+j} f([x_i, x_j], x_1, \ldots, \hat{x}_i, \ldots, \hat{x}_j, \ldots, x_{k+1}).
    \end{align*}
\end{defn}

Let $Z^k(M) = \Ker(d_k) \subset C^k(M)$. Then the compatibility condition $\delta: \mathfrak{g} \to \mathfrak{g} \wedge \mathfrak{g}$ with the Lie bracket can be rewritten as $d \delta = 0$ in the cochain complex with coefficients in the module $M = \mathfrak{g} \wedge \mathfrak{g}$.
\begin{prop}
    If $(\mathfrak{g}, \delta)$ is a Lie bialgebra, then $\delta \in Z^1(\mathfrak{g},\wedge^2\mathfrak{g})$.
\end{prop}

Hence we need to look for a bialgebra structure among 1-cocycles. The easiest way is to construct a coboundary cocycle $\delta = d r$, where $r \in C^0(\mathfrak{g} \wedge \mathfrak{g}) = \mathfrak{g} \wedge \mathfrak{g}$, so that $\delta(x) = x r = (\ad_x \otimes 1 + 1 \otimes \ad_x)(r)$. We will also need an alternative form $\delta(x) = [x_1 + x_2, r]$, where $x_1 = x \otimes 1, x_2 = 1 \otimes x \in \mathcal{U}(\mathfrak{g})^{\otimes 2}$. Thus compatibility of the structures is automatic, but in addition we require that $\delta^*$ be a Lie bracket on $\mathfrak{g}^*$, namely the Jacobi identity should hold. Under what conditions on $r \in \mathfrak{g} \wedge \mathfrak{g}$ does $\delta = d r$ satisfy the Jacobi identity? Necessary and sufficient conditions are called the modified classical Yang--Baxter equation.

Let $r_{12} = r \otimes 1$, $r_{23} = 1 \otimes r$, $r_{13} = P_{23} r_{12} P^{-1}_{23}$, where $P \in \End(\mathfrak{g} \otimes \mathfrak{g}), P(u \otimes v) = v \otimes u$.
\begin{thm}[Proposition 4.20 \cite{STS}]\label{cybe}
    Suppose $[[r, r]] = [r_{12}, r_{13}] + [r_{12}, r_{23}] + [r_{13}, r_{23}] \in (\wedge^3 \mathfrak{g})^{\mathfrak{g}}$. Then $(\mathfrak{g}, \delta = d r)$ is a Lie bialgebra. The converse is also true.
\end{thm}

Here $[[\cdot, \cdot]]: \wedge^k \mathfrak{g} \otimes \wedge^l \mathfrak{g} \to \wedge^{k+l-1} \mathfrak{g}$ is the Schouten bracket, the extension of the Lie bracket on $\mathfrak{g}$ to exterior powers. A special case of this condition is the classical Yang--Baxter equation. The classical Yang--Baxter equation gives a sufficient condition for $r \in \mathfrak{g} \wedge \mathfrak{g}$ to define a coboundary Lie bialgebra:
\begin{equation*}
    [[r, r]] = [r_{12}, r_{13}] + [r_{12}, r_{23}] + [r_{13}, r_{23}] = 0.
\end{equation*}

\subsection{$r$-matrix Poisson brackets}

If $\mathfrak{g}$ is a matrix Lie algebra, for example $\mathfrak{gl}_N$, then the coboundary bialgebra structures from the previous section are called $r$-matrix Poisson brackets, since linear Poisson brackets on $V^*$ and Lie brackets on $V$ correspond bijectively to each other. In fact, it is not necessary that the Lie algebra be matrix-valued (or even finite-dimensional, e.g. a loop algebra), but one often assumes that, first, the Lie algebra is generated by an associative algebra, and second, the matrix representation is natural. In the case of a matrix algebra, the associated linear Poisson bracket is
\begin{equation}\label{r-matrix-linear}
    \{L_1, L_2\}_{lin} = [r_{12}, L_1 + L_2],
\end{equation}
where $t_{ij}$ is a basis of $\mathfrak{gl}^{*}_N$, dual to the basis of matrix units, $L = \sum_{i,j} t_{ij} \otimes E_{ij} \in \mathfrak{gl}^{*}_N \otimes \Mat_N$, $\{L_1, L_2\} \in \mathfrak{gl}^{*}_N \otimes \Mat_N \otimes \Mat_N$, $\{L_1, L_2\}_{i,j,k,l} = \{t_{ij}, t_{kl}\}$. One may think of $L$ as a matrix with coefficients in $\mathfrak{gl}^{*}_N = \Hom(\mathfrak{gl}_N, \mathbb{R}) \subset \Fun(\mathfrak{gl}_N)$, and also as a matrix-valued function on $\mathfrak{gl}_N$, i.e. $L \in \Hom(\mathfrak{gl}_N, \Mat_N)$. This compact way of defining a Poisson bracket on coordinate functions will be used below. By the Leibniz identity, any Poisson bracket is uniquely determined by its values on linear functions, so we often say that a Poisson bracket is defined if its values on coordinate functions are given and the Jacobi identity has been verified.

This $r$-matrix Poisson bracket is the linearization of the Sklyanin bracket
\begin{equation}\label{r-matrix-sklyanin}
    \{L_1, L_2\}_{Sk} = [r_{12}, L_1 L_2].
\end{equation}

In the context of Poisson-Lie group theory, it is more correct to say that the Sklyanin bracket is a bracket on the Poisson-Lie group $GL_N$, and the $r$-matrix bracket is its linearization on $T_e GL_N = \mathfrak{gl}_N$. It is precisely on this pair of brackets that we will focus when constructing $r$-matrix Nambu brackets.

We will also need the definition of compatible Poisson brackets. As with Nambu brackets, we define a Poisson bracket by a bivector field $\pi \in \Gamma(\wedge^2 T M)$, which we call a Poisson bivector, and the bracket is given by $\{f, g\} = \pi(df, dg)$, or more explicitly in coordinates $\{f, g\} = \sum_{i,j} \pi^{ij} \frac{\partial f}{\partial x_i} \frac{\partial g}{\partial x_j}$.
\begin{defn}
    Poisson brackets $\pi_1$ and $\pi_2$ on $M$ are called compatible if for every $\lambda \in \mathbb{R}$, $\pi_1 + \lambda \pi_2$ is a Poisson bracket.
\end{defn}
One of the motivations for studying Nambu brackets is to find families of compatible Poisson brackets.

\begin{prop}
    Let $\pi \in \Gamma(\wedge^3 T M)$ be a Nambu bracket, $f, g \in C^{\infty}(M)$. Then the Poisson brackets $\pi_1 = \iota_{df}\pi$ and $\pi_2 = \iota_{dg}\pi$ are compatible.
\end{prop}
\begin{proof}
    Let $\lambda \in \mathbb{R}$. Then $\pi_1 + \lambda \pi_2 = \iota_{df}\pi + \lambda \iota_{dg}\pi = \iota_{d(f + \lambda g)}\pi$ is also a Poisson bracket by Proposition \ref{substitution}.
\end{proof}

\section{Linear Nambu brackets construction}\label{sec:linear-nambu}

\subsection{Construction of the bracket and basic properties}
Let $V = \mathbb{R}^N$ with a fixed basis $e_1, \ldots, e_N$. Consider $\varphi \in \End(V)^{\wedge 3}$. From now on we identify $\Mat_N = \End(V)$ using the fixed basis in $V$. Define a linear ternary operation on $M = \mathfrak{gl}(V)$

\begin{equation}
    \{t_{i_1 j_1}, t_{i_2 j_2}, t_{i_3 j_3}\} = [\varphi_{123}, L_1 + L_2 + L_3]_{i_1 j_1, i_2 j_2, i_3 j_3},
\end{equation}
where $t_{i j} \in \mathfrak{gl}^{*}(V)$ is the basis dual to the basis of matrix units, $L = (t_{i j})_{i,j = 1, \ldots N} = \sum_{i,j} t_{i j} \otimes E_{i j} \in \mathfrak{gl}^{*}(V) \otimes \End(V), L_i = I \otimes \ldots I \otimes L \otimes I \otimes \ldots \otimes I \in \mathfrak{gl}^{*}(V) \otimes \End(V)^{\otimes l}$, where $L$ appears in the $i$-th tensor factor, $I \in \End(V)$ is the identity matrix, $l$ is a sufficiently large number, clear from the context. In our setting it is almost always 2, 3, or 5. Multiplication in the tensor product is componentwise. The following tensor formula expresses the same operation.

\begin{defn}
    The linear ternary operation defined on coordinate functions by
    \begin{equation}\label{nambu}
    \{L_1, L_2, L_3\} = [\varphi_{123}, L_1 + L_2 + L_3]
    \end{equation}
    is called the $\varphi$-matrix operation, where $\varphi \in \wedge^3 \mathfrak{gl}_N$. It is a Nambu bracket when it satisfies the fundamental identity.
\end{defn}
The operation extends uniquely from linear functions to $C^{\infty}(M)$ by the Leibniz identity:
\begin{equation}
    \{f, g, h\} = \pi^{IJK} \partial_I f \; \partial_J g \; \partial_K h,
\end{equation}
where $I = (i, j), \partial_I = \frac{\partial}{\partial t_{i j}}, \pi^{I J K} = \{t_I, t_J, t_K\}$, and the Einstein summation convention is assumed for repeated indices. We also denote the corresponding trivector field by $\pi = \frac{1}{3!} \pi^{I J K} \partial_{I} \wedge \partial_{J} \wedge \partial_{K} \in \Gamma(\wedge^3 TM)$. The factor accounts for the six permutations of a triple of indices and can be omitted when summing over lexicographically ordered triples $I < J < K$:
\begin{equation}
    \pi = \sum_{I < J < K} \pi^{IJK} \partial_I \wedge \partial_J \wedge \partial_K.
\end{equation}
From now on, whenever an exterior product is present, we will assume summation only over increasing tuples of indices, thus we will not write the sum explicitly and will thus eliminate the symmetric factor.

The skew-symmetry of this operation follows from the skew-symmetry of $\varphi_{123}$:
\begin{equation*}
    \{L_2, L_1, L_3\} = [\varphi_{213}, L_2 + L_1 + L_3] = [-\varphi_{123}, L_1 + L_2 + L_3] = -\{L_1, L_2, L_3\},
\end{equation*}
where $\varphi_{213} = P_{12} \varphi_{123} P_{12}, P = \sum_{i, j} E_{i j} \otimes E_{j i} \in \End(V)^{\otimes 2}$ is the permutation operator
\begin{equation*}
    P(u \otimes v) = v \otimes u.
\end{equation*}

It is clear that the conditions on the skew-symmetry of $\varphi$ is only sufficient and can be dropped. The necessary condition is that the differences $\varphi_{123} - (-1)^{\sigma} \varphi_{\sigma(1) \sigma(2) \sigma(3)}$ belong to the invariants of the action $(\End(V \otimes V \otimes V))^{\mathfrak{g}}$ for $\sigma \in S_3$. The action of $\mathfrak{g} = \mathfrak{gl}(V)$ is the standard one under the identification $\End(V) = \mathfrak{g}$, namely
\begin{equation}
\begin{array}{c}
    x (u \otimes v \otimes w) = (\ad_x \otimes \id \otimes \id + \id \otimes \ad_x \otimes \id + \id \otimes \id \otimes \ad_x) (u \otimes v \otimes w) \\= [x, u] \otimes v \otimes w + u \otimes [x, v] \otimes w + u \otimes v \otimes [x, w].
\end{array}
\end{equation}

From now on we assume $\varphi$ is skew-symmetric unless stated otherwise.

\subsection{Study of the $r$-matrix form of the bracket}

We specialize to the class of linear brackets for $\varphi_{123} = \varphi(r) =  r_{12} - r_{13} + r_{23} = A(r \otimes I)$, where $r \in \End(V)^{\wedge 2}$, $A: \End(V)^{\otimes 3} \to \End(V)^{\wedge 3}$ is the antisymmetrizer.

\begin{defn}
    The linear ternary operation defined on coordinate functions by
    \begin{equation}\label{r-matrx-nambu}
    \{L_1, L_2, L_3\} = [r_{12} - r_{13} + r_{23}, L_1 + L_2 + L_3]
    \end{equation}
    is called the $r$-matrix operation, where $r \in \wedge^2 \mathfrak{gl}_N$.
\end{defn}

We call it an $r$-matrix Nambu bracket when it satisfies the fundamental identity \eqref{fund_id_3}. This choice of $\varphi$ is motivated by the fact that $R \otimes I$ is a (trivial) solution of the Zamolodchikov tetrahedron equation (for more details on the tetrahedron equation, see \cite{Talalaev2021}),
\begin{equation*}
    \Phi_{123} \Phi_{145} \Phi_{246} \Phi_{356} = \Phi_{356} \Phi_{246} \Phi_{145} \Phi_{123},
\end{equation*}
if $R$ is a solution of the Yang--Baxter equation. At the same time, the classical Yang--Baxter equation is the linearization of the Yang--Baxter equation, namely let $R = \Id + \hbar r + o(\hbar^2)$. The Yang--Baxter equation is an equation for $R \in \End(V)^{\otimes 2}$:
\begin{equation*}
    R_{12} R_{13} R_{23} = R_{23} R_{13} R_{12}.
\end{equation*}
Expand this equation with respect to the formal parameter $\hbar$. Then for $\hbar^0$ and $\hbar^1$ we get identities, and the first nontrivial equality appears at $\hbar^2$, and it is called the classical Yang--Baxter equation:
\begin{equation*}
    [r_{12}, r_{13}] + [r_{12}, r_{23}] + [r_{13}, r_{23}] = 0.
\end{equation*}

For a skew-symmetric $r$, let $\{\cdot,\cdot\}_{lin}$ and
$\{\cdot,\cdot\}_{Sk}$ denote the binary operations defined by
\eqref{r-matrix-linear} and \eqref{r-matrix-sklyanin}, respectively.
Either operation may fail the Jacobi identity for a general $r$.
Both are Poisson when $r$ satisfies the modified classical
Yang--Baxter equation, i.e., when $\operatorname{CYB}(r)$ is
invariant under the adjoint action.

Our basic example is the classical Drinfeld--Jimbo $r$-matrix, more precisely its skew-symmetric part $r_{12} - r_{21}$,
\begin{equation*}
    r_{DJ} = \sum_{i < j} E_{i j} \wedge E_{j i}.
\end{equation*}
Fixing a trace function in an $r$-matrix Nambu bracket
\eqref{r-matrx-nambu} yields a family of compatible Poisson brackets.
The formulas for the first two members involve the linear and Sklyanin
bivectors defined above.

\begin{thm}\label{r-matrix-brackets}
    Let $\{\cdot, \cdot, \cdot\}$ be an $r$-matrix Nambu bracket. Then we have a family of compatible Poisson brackets
    $\{f, g\}_n = \{f, g, \Tr L^n\}$, where
    \begin{align*}
      \{L_1,L_2\}_1
      & =N\{L_1,L_2\}_{lin}
      -[\Tr_3(r_{13}),L_1]+[\Tr_3(r_{23}),L_2]\\
      & =N\{L_1,L_2\}_{lin}
      + I\wedge[\Tr_2(r_{12}),L_1],\\
      \{L_1,L_2\}_2
      & =(2\Tr L)\{L_1,L_2\}_{lin}
      +2I\wedge\{L,\Tr L\}_{Sk}.
    \end{align*}
\end{thm}
\begin{proof}
    Write $r_{12} = \sum_{\alpha} A_{\alpha} \otimes B_{\alpha} = A_{\alpha} \otimes B_{\alpha}$, where repeated indices are understood as summed, and let $N = \dim V$. We will also use the cyclic property of the trace whenever the operators act only in the space over which the trace is taken, namely $\Tr_{k} (X_k Y_k) = \Tr_{k} (Y_k X_k)$. For example, it follows that $\Tr_k [X_k, Y_k] = 0$. Let $I \in \End(V)$ denote the identity matrix.
    We will prove the statement by direct computation using matrix techniques.
    \begin{align*}
        &\{L_1, L_2\}_1 = \{L_1, L_2, \Tr L\} = \Tr_3 (\{L_1, L_2, L_3\}) = \\
        &= \Tr_3 [r_{12} - r_{13} + r_{23}, L_1 + L_2 + L_3] = S_{12} - S_{13} + S_{23},
    \end{align*}
    where $S_{ij} = \Tr_3 [r_{ij}, L_1 + L_2 + L_3]$. Let us compute each term separately.
    \begin{align*}
        S_{12} = \Tr_3 [r_{12}, L_1 + L_2] + \Tr_3 [r_{12}, L_3] = N \{L_1, L_2\} + \Tr_3 ([r_{12}, L_3]) = N \{L_1, L_2\},
    \end{align*}
    where $[r_{12}, L_3] = 0$ because these matrices act in different copies of the space $V$.
    \begin{align*}
        S_{13} = \Tr_3[r_{13}, L_1 + L_2 + L_3] = \Tr_3 [r_{13}, L_1 + L_3] = \Tr_3 [r_{13}, L_1] + \Tr_3 [r_{13}, L_3].
    \end{align*}
    Compute $\Tr_3 [r_{13}, L_3]$.
    \begin{align*}
        \Tr_3 [r_{13}, L_3] = \Tr_3 (A_{\alpha} \otimes I \otimes B_{\alpha} L - A_{\alpha} \otimes I \otimes L B_{\alpha}) = \Tr_3 (A_{\alpha} \otimes I \otimes [B_{\alpha}, L]) = 0.
    \end{align*}
    Compute $\Tr_3 [r_{13}, L_1]$.
    \begin{align*}
        \Tr_3 [r_{13}, L_1] = \Tr_3 ([A_{\alpha}, L] \otimes I \otimes B_{\alpha}) = (\Tr B_{\alpha}) [A_{\alpha}, L] \otimes I = [\Tr_{3} r_{13}, L_1].
    \end{align*}
    Hence $S_{13} = [\Tr_3 r_{13}, L_1]$. The computation of $S_{23}$ is analogous up to a relabeling of the space index, so $S_{23} = [\Tr_3 r_{23}, L_2]$. Thus we have computed $\{L_1, L_2\}_1$.

    Now compute $\{L_1, L_2\}_2$.
    \begin{align*}
        \{L_1, L_2\}_2 &= \Tr_3 \{L_1, L_2, L^2_3\} = \Tr_3 (L_3 \{L_1, L_2, L_3\} + \{L_1, L_2, L_3\} L_3) =  K_{12} - K_{13} + K_{23}, 
    \end{align*}
    where $K_{ij} = \Tr_3 (L_3 [r_{ij}, L_i + L_j] + [r_{ij}, L_i + L_j] L_3)$.
    It is simplest to compute $K_{12} = 2 (\Tr L) \{L_1, L_2\}_{lin}$.
    Now compute $K_{13}$,
    \begin{align*}
        K_{13} &= \Tr_3 (L_3 [r_{13}, L_1 + L_3] + [r_{13}, L_1 + L_3] L_3) = \\ &= \Tr_3 (L_3 r_{13} L_1 - L_3 L_1 r_{13} + \cancel{L_3 r_{13} L_3} - \cancel{L^2_3 r_{13}} + r_{13} L_1 L_3 - L_1 r_{13} L_3 + \cancel{r_{13} L^2_3} - \cancel{L_3 r_{13} L_3}) = \\ &= \Tr_3 (r_{13} L_1 L_3 - L_1 L_3 r_{13} + r_{13} L_1 L_3 - L_1 L_3 r_{13}) = 2 \Tr_3 [r_{13}, L_1 L_3] = 2 \Tr_3 \{L_1, L_3\}_{Sk},
    \end{align*}
    where we used the cyclic property of the trace, and also $L_1 L_3 = L_3 L_1$, because the operators lie in different copies of the space, while the $t_{ij}$ are functions and therefore commute.
    The computation of $K_{23}$ is completely analogous, giving $K_{23} = 2 \Tr_3 \{L_2, L_3\}_{Sk} = 2 \{L_2, \Tr_3 L_3\}_{Sk}$.
    Therefore
    \begin{align*}
        \{L_1, L_2\}_2 &= 2 (\Tr L) \{L_1, L_2\}_{lin} - 2 \{L_1, \Tr_3(L_3)\}_{Sk} + 2 \{L_2, \Tr_3(L_3)\}_{Sk} = \\ &= 2 (\Tr L) \{L_1, L_2\}_{lin} + 2 I \wedge \{L, \Tr L\}_{Sk}.
    \end{align*}
\end{proof}

\begin{cor}\label{cor-bracket}
    Suppose the assumptions of Theorem \ref{r-matrix-brackets} hold. Then the following refinements are true:
    \begin{itemize}
        \item If $\Tr_2(r_{12}) = 0$, then
        \begin{equation*}
            \{f, g\}_1 = N \{f, g\}_{lin}.
        \end{equation*}
        \item If $N = \dim V = 2$, then \begin{equation*}
        \{L_1, L_2, \det L\} = (\Tr L) \{L_1, L_2\}_{lin} + I \wedge ((\Tr L) \{L_1, \Tr L\}_{lin} - \{L_1, \Tr L\}_{Sk}).
        \end{equation*}
    \end{itemize}
\end{cor}
\begin{proof}
The first corollary is obvious, since the second term in the expression for $\{L_1, L_2\}_1$ vanishes. For the second corollary, write
\begin{equation*}
    \det L = \frac{(\Tr L)^2 - \Tr L^2}{2}.
\end{equation*}
Then
\begin{align*}
    \{L_1, L_2, \det L\} &= \frac{1}{2} (2 (\Tr L) \{L_1, L_2\}_1 - \{L_1, L_2\}_2) = \\ &= (\Tr L) \{L_1, L_2\}_{lin} + I \wedge ((\Tr L) \{L_1, \Tr L\}_{lin} - \{L_1, \Tr L\}_{Sk}).
\end{align*}
\end{proof}

\begin{cor}\label{weak-FI}
    Let $r \in \mathfrak{sl}(V)^{\wedge 2}$ define an $r$-matrix Nambu bracket on $\mathfrak{gl}(V)$. Then $r$ satisfies the modified classical Yang--Baxter equation, namely $[[r, r]] \in (\mathfrak{g}^{\wedge 3})^{\mathfrak{g}}$, where $\mathfrak{g} = \mathfrak{gl}(V)$.
\end{cor}
\begin{proof}
    Since $r \in \mathfrak{sl}(V)^{\wedge 2}$, we have $\Tr_2(r_{12}) = 0$. Therefore $\{L_1, L_2\}_1 = N \{L_1, L_2\}_{lin} = N [r_{12}, L_1 + L_2]$ is a Poisson bracket by Corollary \ref{cor-bracket}. By Theorem \ref{cybe} we know that the Jacobi identity for this Poisson bracket is equivalent to the modified classical Yang--Baxter equation.
\end{proof}

In fact, we can weaken the conditions on the Nambu bracket if we only want to obtain a linear Poisson bracket.

\begin{prop}\label{not-deco}
    Let $M = \mathfrak{gl}_N$ and let $\pi \in \Gamma(\wedge^3 T M)$ be linear (in fixed basis) and define a 3-Lie algebra structure on $M^{*}$, but not necessarily a Nambu bracket structure (see the counterexample in Proposition \ref{counterexample}). Then $\{f,g\}_1=\{f,g,\Tr L\}$ is a Poisson bracket.
\end{prop}
\begin{proof}
    We have that the fundamental identity holds for linear functions. $\Tr L$ is also a linear function. As is well known, to verify the Poisson identity it is enough to do so on linear functions, therefore $\{\cdot, \cdot\}_1$ is a Poisson bracket.
\end{proof}

\subsection{Linear Nambu brackets on $\mathfrak{gl}_2$}

Let us introduce the standard notation for the generators of $\mathfrak{sl}_2 \subset \mathfrak{gl}_2$, namely $e = E_{12}, h = E_{11} - E_{22}, f = E_{21}$.

Our current task is to determine for which $r \in \wedge^2 \mathfrak{gl}_2$ the expression $[\varphi_{123}, L_1 + L_2 + L_3]$ satisfies the fundamental identity, where $\varphi_{123} = \varphi(r) = r_{12} - r_{13} + r_{23}$. In this section $\mathfrak{g} = \mathfrak{gl}_2$. Moreover, in explicit computations we use the notation $\varphi_{ijk}$ for the components of the tensor $\varphi \in \wedge^3 \mathfrak{gl}_2$, not for the labels of spaces. Recall that brackets of the form $[\varphi, L_1 + L_2 + L_3]$ are called $\varphi$-matrix, and their specialization $\varphi = r_{12} - r_{13} + r_{23}$ are called $r$-matrix brackets.

\begin{thm}\label{thm-gl2}
    Let $\varphi \in \wedge^3 \mathfrak{gl}_2$. Then $[\varphi, L_1 + L_2 + L_3]$ defines a Nambu bracket on $\mathfrak{gl}_2$.
\end{thm}
\begin{proof}
    Write $(t_1,t_2,t_3,t_4)=(t_{11},t_{12},t_{21},t_{22})$ and let $\Omega=d t_1\wedge d t_2\wedge d t_3\wedge d t_4$. Computing the commutator on each of the three factors of a basis wedge gives
    \begin{align}
        \{t_{1}, t_{2}, t_{3}\} &= \varphi_{134} t_{2} + \varphi_{124} t_{3}, \\
        \{t_{1}, t_{2}, t_{4}\} &= -\varphi_{124} t_{1} + (\varphi_{123} + \varphi_{234}) t_{2} + \varphi_{124} t_{4}, \\
        \{t_{1}, t_{3}, t_{4}\} &= \varphi_{134} t_{1} + (\varphi_{123} + \varphi_{234}) t_{3} - \varphi_{134} t_{4}, \\
        \{t_{2}, t_{3}, t_{4}\} &= \varphi_{134} t_{2} + \varphi_{124} t_{3}.
    \end{align}

    Put $b=\varphi_{124}$, $c=\varphi_{134}$, $s=\varphi_{123}+\varphi_{234}$, $p=t_1-t_4$, and $u=ct_2+bt_3$. For $\eta=t_3\,d t_2-t_2\,d t_3$, the dual 1-form and its derivative are
    \begin{align*}
        \omega=\iota_\pi\Omega &= -u\,d p+p\,d u+s\eta,\\
        d\omega &= 2\,d p\wedge d u-2s\,d t_2\wedge d t_3.
    \end{align*}
    Since $\eta\wedge d u=u\,d t_2\wedge d t_3$, the two mixed terms in $\omega\wedge d\omega$ cancel and all other terms vanish. Corollary \ref{corollary_fund} now gives the fundamental identity.
\end{proof}

Note that the invariants of the action of $\mathfrak{g}$ on $\wedge^3 \mathfrak{g}$ give the zero Nambu bracket. Therefore, in fact, the space of distinct $\varphi$-matrix Nambu brackets has dimension $\dim \wedge^3 \mathfrak{g} - \dim (\wedge^3 \mathfrak{g})^{\mathfrak{g}} = \binom{4}{3} - 1 = 3$. We also note that the space of $\varphi$-matrix Nambu brackets is linear.

It follows that all $r$-matrix Poisson brackets on $\mathfrak{g}$ with the condition $\Tr_2 (r_{12}) = 0$ are compatible.

\begin{cor}
    Let $r \in \wedge^2 \mathfrak{g}$ and $\Tr_2(r_{12}) = 0$. Then $[r_{12}, L_1 + L_2]$ defines a Poisson bracket, and these brackets (for different choice of $r$) are compatible with each other.
\end{cor}

\begin{proof}
    Construct $\varphi = \varphi(r)$ and $\tilde{\varphi} = \varphi(\tilde{r})$. Then $\varphi_{\lambda} = \varphi + \lambda \tilde{\varphi}$ will define a Nambu bracket for any value of the parameter $\lambda$. Construct the Nambu bracket for $\varphi_{\lambda}$, substitute $\Tr L$, and by Theorem \ref{r-matrix-brackets} obtain a linear Poisson bracket for $r + \lambda \tilde{r}$. Therefore the Poisson brackets are compatible.
\end{proof}

This simple fact also admits an independent proof that does not use Nambu brackets.
\begin{proof}
    Let $r \in \wedge^2 \mathfrak{gl}_N$ and $\Tr_2 (r_{12}) = 0$. We show that $r \in \wedge^2 \mathfrak{sl}_N$. It is clear that $\mathfrak{gl}_N = \mathfrak{sl}_N \oplus \mathbb{R} I$, where $I$ is the identity matrix. Then $\wedge^2 \mathfrak{gl}_N = \wedge^2 \mathfrak{sl}_N + \mathbb{R} I \wedge \mathfrak{sl}_N$. Thus $r = r_{0} + I \wedge a_0$, where $r_0 \in \wedge^2 \mathfrak{sl}_N, a_0 \in \mathfrak{sl}_N$. Compute $\Tr_2(r_{12}) = \Tr_2(r_0) + \Tr_2(I \otimes a_0) - \Tr_2(a_0 \otimes I) = 0 + 0 - N a_0 = - N a_0$. Since by assumption $\Tr_2(r_{12}) = 0$, we get $a_0 = 0$, hence $r = r_0 \in \wedge^2 \mathfrak{sl}_N$.

    It is known that to verify the Jacobi identity, it suffices to check that $[[r,r]] \in (\wedge^3 \mathfrak{g})^{\mathfrak{g}}$. In our case $r \in \wedge^2 \mathfrak{sl}_2$, hence $[[r,r]] \in \wedge^3 \mathfrak{sl}_2 = \mathbb{R} e \wedge h \wedge f$. Thus $[[r,r]] = \lambda e \wedge h \wedge f$ for some $\lambda \in \mathbb{R}$. But this element is obviously annihilated by the action of $\mathfrak{g} = \mathfrak{gl}_2$, so the modified Yang--Baxter equation is automatic.

    Note that the condition $\Tr_{2}(r_{12}) = 0$ is not equivalent to $r \in \otimes^2 \mathfrak{sl}_N$ when $r$ is not skew-symmetric (i.e. $r_{12} \ne -r_{21}$), for example, $r = I \otimes E_{12}$. This highlights the specificity of considering only skew-symmetric classical $r$-matrices.
\end{proof}

Nevertheless, this fact is not totally trivial, because $[[r, r]] \ne 0$ for some $r \in \wedge^2 \mathfrak{sl}_2$, for example, if $r = e \wedge f$, then $[[r, r]] = h \wedge e \wedge f \ne 0$.

We now show that all $\varphi$-matrix Nambu brackets on $\mathfrak{gl}_2$ are in fact $r$-matrix ones.

\begin{prop}
    All $\varphi$-matrix Nambu brackets on $\mathfrak{gl}_2$ are $r$-matrix ones.
\end{prop}

\begin{proof}
    We have already understood that both spaces are linear and have the same dimension. Therefore it is enough to provide three $r$-matrices giving linearly independent Nambu brackets. Take $r_1 = E_{11} \wedge E_{12}$, then $\varphi_1 = \varphi(r_1) = E_{11} \wedge E_{12} \wedge E_{22}$, which corresponds to the coefficients $\varphi_{123} = \varphi_{134} = \varphi_{234} = 0, \varphi_{124} = 1$ in the formulas for the bracket from the proof of Theorem \ref{thm-gl2}. Let $r_2 = E_{11} \wedge E_{21}$. Then $\varphi_2 = \varphi(r_2) = E_{11} \wedge E_{21} \wedge E_{22}$. Hence the bracket $\varphi_2$ corresponds to the coefficients $\varphi_{123} = \varphi_{124} = \varphi_{234} = 0, \varphi_{134} = 1$. Finally, take $r_3 = E_{12} \wedge E_{21}$, then $\varphi_3 = \varphi(r_3) = E_{11} \wedge E_{12} \wedge E_{21} + E_{12} \wedge E_{21} \wedge E_{22}$, i.e. $\varphi_{123} = \varphi_{234} = 1, \varphi_{124} = \varphi_{134} = 0$. It is clear that the three resulting $r$-matrix brackets are linearly independent.
\end{proof}

\subsection{Main example}

Let $r_{DJ} = E_{12} \wedge E_{21}$. The corresponding $r$-matrix Nambu bracket is given by
\begin{equation*}
    \pi = 2 t_{12} \partial_{11} \wedge \partial_{12} \wedge \partial_{22} + 2 t_{21} \partial_{11} \wedge \partial_{21} \wedge \partial_{22} = 2 \partial_{11} \wedge (t_{12} \partial_{12} + t_{21} \partial_{21}) \wedge \partial_{22}.
\end{equation*}
Hence $\{t_{11}, t_{12}, t_{22}\} = 2 t_{12}, \{t_{11}, t_{21}, t_{22}\} = 2 t_{21}$.
Let us denote the partial Euler vector field $v = t_{12} \partial_{12} + t_{21} \partial_{21}$. Then
\begin{equation*}
    \pi = 2 \partial_{11} \wedge v \wedge \partial_{22}.
\end{equation*}
Substituting $f = \Tr L = t_{11} + t_{22}$ gives
\begin{equation*}
    \iota_{df}\pi = 2 \partial_{11} \wedge v - 2 \partial_{22} \wedge v = 2 (\partial_{11} - \partial_{22}) \wedge v.
\end{equation*}
We obtain $\iota_{df}\pi = 2 \{\cdot, \cdot\}_{lin}$, which agrees with Theorem \ref{r-matrix-brackets}. Also, substituting $g = \det L = t_{11} t_{22} - t_{12} t_{21}$ yields $\iota_{dg}\pi = 2 \{ \cdot, \cdot \}_{Sk}$. Thus, using the general Nambu bracket, we have shown that the linear $r$-matrix Poisson bracket and the Sklyanin bracket are compatible. This is a well-known fact in arbitrary dimension (Theorem 13 \cite{STS83}). But in $\mathfrak{gl}_2$ they are linked even more rigidly.

\subsection{Restriction to $\mathfrak{sl}_2$}

We show that only the zero $\varphi$-matrix Nambu bracket restricts to $\mathfrak{sl}_2$.

\begin{lemma}\label{sl-restriction}
    Let $\pi$ be a $\varphi$-matrix Nambu bracket on $\mathfrak{gl}_2$ for which the restriction to $\mathfrak{sl}_2$ is well-defined. Then $\pi = 0$.
\end{lemma}

\begin{proof}
    Put $q= \Tr L = t_{11}+t_{22}$. Since the ideal of functions vanishing on $\mathfrak{sl}_2$ is generated by $q$, the bracket descends if and only if $\{q,f,g\}|_{q=0}=0$ for all $f,g$. The coordinate formulas above give
    \begin{align*}
        \{q,t_{12},t_{21}\}|_{q=0} &= 2\varphi_{134}t_{12}+2\varphi_{124}t_{21},\\
        \{q,t_{11},t_{12}\}|_{q=0} &= -2\varphi_{124}t_{11}+(\varphi_{123}+\varphi_{234})t_{12}.
    \end{align*}
    Since $t_{11},t_{12},t_{21}$ are independent on $q=0$, vanishing forces $\varphi_{134}=\varphi_{124}=\varphi_{123}+\varphi_{234}=0$. The coordinate formulas then give $\pi=0$ on $\mathfrak{gl}_2$.
\end{proof}

\section{Quadratic and cubic Nambu brackets construction}\label{cubic-section}

\subsection{Definition and basic properties}

In the standard theory of Poisson-Lie groups, the $r$-matrix Poisson bracket is the linearization of the Sklyanin bracket on the Poisson-Lie group. In this connection, it is natural to consider a cubic operation whose linearization yields the $\varphi$-matrix linear Nambu bracket.

\begin{defn}
    Let $\varphi \in \wedge^3 \mathfrak{gl}_N$. Then we call
    \begin{equation}\label{cubic-trivector}
        \{L_1, L_2, L_3\} = [\varphi, L_1 L_2 L_3]
    \end{equation}
    the cubic ternary operation.
\end{defn}

\begin{defn}
    Let $\varphi \in \wedge^3 \mathfrak{gl}_N$. Then we call
    \begin{equation*}
        \{L_1, L_2, L_3\} = [\varphi, L_1 L_2 + L_1 L_3 + L_2 L_3]
    \end{equation*}
    the quadratic ternary operation.
\end{defn}

Indeed, if $L = I + \hbar K$, where $\hbar$ is a formal parameter, then
\begin{equation*}
    [\varphi, L_1 L_2 L_3] = \hbar [\varphi, K_1 + K_2 + K_3] + \hbar^2 [\varphi, K_1 K_2  + K_1 K_3 + K_2 K_3] + \hbar^3 [\varphi, K_1 K_2 K_3].
\end{equation*}

\begin{prop}\label{cubic_bracket}
Let $\varphi \in \wedge^3 \mathfrak{gl}_N$, and suppose the cubic operation satisfies the fundamental identity. Then its linear part also satisfies the fundamental identity.
\end{prop}
\begin{proof}
    In fixed matrix coordinates, evaluate the cubic tensor at $L=I+\hbar K$. Its coefficients start with $\hbar\pi_{\rm lin}(K)$, and their first derivatives with respect to the entries of $L$ start with the constant first derivatives of $\pi_{\rm lin}$. Thus the coefficient of $\hbar$ in the differential identity of Theorem \ref{thm-vaisman} is exactly that identity for $\pi_{\rm lin}$.

    The algebraic Pl\"ucker relations are homogeneous quadratic equations in the tensor coefficients. The coefficient of $\hbar^2$ in each relation is therefore the corresponding relation for $\pi_{\rm lin}$. The two conditions of Theorem \ref{thm-vaisman} now give the full fundamental identity for $\pi_{\rm lin}$.
\end{proof}

It is well-known that Poisson bracket on Poisson-Lie group is a multiplicative bivector field.

\begin{defn}
    Polyvector field $\pi \in \Gamma(\wedge^n TG)$ where $G$ is a Lie group is called multiplicative if
    $$
    \pi(g h) = (L_{g})_{*} \pi(h) + (R_{h})_{*} \pi(g), 
    $$
    where $L_g: G \to G$ and $R_g: G \to G$ are left and right multiplication.
\end{defn}
This definition is used to generalize Poisson-Lie groups to Nambu-Lie groups.

\begin{defn}[\cite{VaismanNambuLie}]
    Nambu-Lie group $G$ is a Lie group with a multiplicative Nambu bracket $\pi \in \Gamma(\wedge^3 T G)$
\end{defn}

\begin{prop}\label{prop:coboundary-mult}
    For any $\varphi \in \wedge^3 \mathfrak{gl}_N$, the cubic trivector \eqref{cubic-trivector} on $GL_N(\mathbb{R})$ is multiplicative. If it satisfies the fundamental identity, it makes $GL_N(\mathbb{R})$ a Nambu-Lie group.
\end{prop}

\begin{proof}
Multiplicativity follows from the commutator product rule.
$$
[\varphi, (L M)_1 (L M)_2 (L M)_3] = L_1 L_2 L_3 [\varphi, M_1 M_2 M_3] + [\varphi, L_1 L_2 L_3] M_1 M_2 M_3,
$$
which is true, because $(L M)_1 (L M)_2 (L M)_3 = L_1 L_2 L_3 M_1 M_2 M_3$ and the commutator satisfies the Leibniz identity.    
\end{proof}

There is a following restrictive condition on multiplicative Nambu trivectors.
\begin{thm}\label{thm:gl-vanishing}
For $N\ge3$, every multiplicative Nambu trivector on $GL_N(\mathbb R)$ vanishes identically.
\end{thm}
\begin{proof}
By the structure theorem of Grabowski and Marmo \cite{GrabowskiMarmo}, a nonzero multiplicative Nambu tensor of order three on a connected Lie group requires an ideal of dimension two, three or four in its Lie algebra. The only ideals of $\mathfrak{gl}_N(\mathbb R)=\mathbb RI\oplus\mathfrak{sl}_N(\mathbb R)$ are zero, the centre, $\mathfrak{sl}_N$ and the whole algebra. For $N\ge3$ their dimensions are $0,1,N^2-1,N^2$. Hence the tensor vanishes on the identity component $G^0$.

To cover the other connected component of the group, let $c(g)=(R_{g^{-1}})_*\pi(g)$ be the trivialization of the trivector field $\pi$, where $c: GL_N \xrightarrow{} \wedge^3 \mathfrak{gl}_N $. Multiplicativity gives $c(gh)=c(g)+\operatorname{Ad}_g c(h)$. For $h\in G^0$ we have $c(h) = 0$. Computing $c(hg)=c(g(g^{-1}hg))$ in two ways gives $\operatorname{Ad}_h c(g)=c(g)$. If $c(g)$ were nonzero, its decomposability would give a three-dimensional support $E \subset \mathfrak{gl}_N$ (the span of its three wedge factors). Now we show that $E$ is $G^0$-invariant. If $X \in E$ then $X \wedge c(g) = 0$. So $0 = \operatorname{Ad}_h(X \wedge c(g)) = (\operatorname{Ad}_h X) \wedge (\operatorname{Ad}_h c(g)) = (\operatorname{Ad}_h X) \wedge c(g)$. Finally, we have that $\operatorname{Ad}_h X \in E$. Hence we obtain a three-dimensional ideal, which is impossible for $N \ge 3$.
\end{proof}

In particular, the cubic trivector \eqref{cubic-trivector} of Proposition~\ref{prop:coboundary-mult} satisfies the full fundamental identity for $N\ge3$ only when it vanishes.

\begin{defn}
    Nambu brackets $\pi_1, \pi_2 \in \Gamma(\wedge^n T M)$ are called compatible if for every parameter $\lambda \in \mathbb{R}$ the sum $\pi_1 + \lambda \pi_2$ is also a Nambu bracket.
\end{defn}

Note that if Nambu brackets of order 3, $\pi_1$ and $\pi_2$, are compatible, then for any function $f \in C^{\infty}(M)$ the Poisson brackets $\iota_{df}\pi_1$ and $\iota_{df}\pi_2$ are compatible. Thus we obtain another tool for constructing compatible Poisson brackets.

\begin{prop}
Let $\varphi \in \wedge^3 \mathfrak{gl}_N$. If the quadratic $\varphi$-matrix operation satisfies the fundamental identity, then
\begin{itemize}
    \item the linear $\varphi$-matrix Nambu bracket also satisfies the fundamental identity,
    \item the quadratic and linear $\varphi$-matrix Nambu brackets on $\mathfrak{gl}_N$ are compatible.
\end{itemize}
\end{prop}
\begin{proof}
    The first statement is proved similarly to Proposition \ref{cubic_bracket} by expanding in the parameter $\hbar$. We obtain in the first order in $\hbar$
    \begin{equation*}
        2 [\varphi, K_1 + K_2 + K_3] = 2 \{K_1, K_2, K_3\}_{lin}.
    \end{equation*}
    The proof of the second statement is analogous to the compatibility of the linear $r$-matrix Poisson bracket and the Sklyanin bracket \cite{STS83}. Consider the affine automorphism $\mathfrak{gl}_N \to \mathfrak{gl}_N$ given by translation by the constant $\lambda I \in \mathfrak{gl}_N$, and denote by $\{L_1, L_2, L_3\}_{\lambda}$ the quadratic Nambu bracket in the new coordinates, where $\tilde{L} = L + \lambda I$. It still satisfies the fundamental identity. At the same time,
    \begin{equation*}
        \{L_1, L_2, L_3\}_{\lambda} = \{L_1, L_2, L_3\} + 2 \lambda \{L_1, L_2, L_3\}_{lin}.
    \end{equation*}
    The left-hand side is a Nambu bracket for every $\lambda$, hence the brackets on the right are compatible.
\end{proof}

\subsection{In search of the analogue of CYBE}

In this section we rewrite the fundamental identity on linear functions for the linear and cubic operations. It already gives us 3-Lie algebra structure on the dual space (given linear $\varphi$-matrix operation). The full fundamental identity also requires the decomposability condition of Theorem \ref{thm-vaisman} which we don't discuss in this Section. Here $\varphi_{ijk}$ denotes the indices of the tensor factors.

    Expand the left-hand side of the (linear) fundamental identity for $\{L_1, L_2, L_3\} = [\varphi_{123}, L_1 L_2 L_3]$. First, use the Leibniz identity and the fact that $\varphi$ is a constant tensor:
    \begin{align*}
        \{L_1, L_2, \{L_3, L_4, L_5\}\} = \{L_1, L_2, [\varphi_{345}, L_3 L_4 L_5]\} = [\varphi_{345}, \{L_1, L_2, L_3 L_4 L_5\}].
    \end{align*}
    Now compute $\{L_1, L_2, L_3 L_4 L_5\}$ using the Leibniz identity:
    \begin{equation*}
        \{L_1, L_2, L_3 L_4 L_5\} = [\varphi_{123}, L_1 L_2 L_3] L_4 L_5 + L_3 [\varphi_{124}, L_1 L_2 L_4] L_5 + L_3 L_4 [\varphi_{125}, L_1 L_2 L_5].
    \end{equation*}
    It remains to use the fact that operators in different tensor factors commute, hence $[\varphi_{123}, L_1 L_2 L_3] L_4 L_5 = [\varphi_{123}, \mathcal{L}]$, where we denote $\mathcal{L} = L_1 L_2 L_3 L_4 L_5$. Hence
    \begin{equation*}
        \{L_1, L_2, L_3 L_4 L_5\} = [\varphi_{123} + \varphi_{124} + \varphi_{125}, \mathcal{L}].
    \end{equation*}
    Therefore
    \begin{equation*}
        \LHS = \{L_1, L_2, \{L_3, L_4, L_5\}\} = [\varphi_{345}, [\varphi_{123} + \varphi_{124} + \varphi_{125}, \mathcal{L}]].
    \end{equation*}
    Rewrite the right-hand side of the fundamental identity \eqref{fund_id_3}, reordering the arguments for uniformity:
    \begin{equation*}
        \RHS = \{\{L_1, L_2, L_3\}, L_4, L_5\} + \{\{L_1, L_2, L_4\}, L_5, L_3\} + \{\{L_1, L_2, L_5\}, L_3, L_4\}.
    \end{equation*}
    Now compute one of the terms in the right-hand side, the remaining ones are obtained by cyclic permutation of the indices $(345)$. For example, compute $\{\{L_1, L_2, L_3\},L_4, L_5\}$.
    \begin{equation*}
        \{\{L_1, L_2, L_3\},L_4, L_5\} = [\varphi_{123}, \{L_1 L_2 L_3, L_4, L_5\}] = [\varphi_{123}, [\varphi_{145} + \varphi_{245} + \varphi_{345}, \mathcal{L}]].
    \end{equation*}
    Hence we can write the whole right-hand side:
    \begin{align*}
        \RHS &= [\varphi_{123}, [\varphi_{145} + \varphi_{245} + \varphi_{345}, \mathcal{L}]] + [\varphi_{124}, [\varphi_{153} + \varphi_{253} + \varphi_{453}, \mathcal{L}]] + [\varphi_{125}, [\varphi_{134} + \varphi_{234} + \varphi_{534}, \mathcal{L}]] \\ &\quad= [\varphi_{123}, [\varphi_{145} + \varphi_{245} + \varphi_{345}, \mathcal{L}]] + [\varphi_{124}, [- \varphi_{135} - \varphi_{235} + \varphi_{345}, \mathcal{L}]] + [\varphi_{125}, [\varphi_{134} + \varphi_{234} + \varphi_{345}, \mathcal{L}]].
    \end{align*}
    This equation can be simplified further by grouping together all terms containing $\varphi_{345}$ on the left-hand side:
    \begin{align*}
    \widetilde{\LHS} &= [\varphi_{345}, [\varphi_{123} + \varphi_{124} + \varphi_{125}, \mathcal{L}]] - [\varphi_{123} + \varphi_{124} + \varphi_{125}, [\varphi_{345}, \mathcal{L}]] = \\ &= -[\mathcal{L}, [\varphi_{345}, \varphi_{123} + \varphi_{124} + \varphi_{125}]].
    \end{align*}
    Then rewrite the right-hand side:
    \begin{align*}
        \widetilde{RHS} = [\varphi_{123}, [\varphi_{145} + \varphi_{245}, \mathcal{L}]] - [\varphi_{124}, [\varphi_{135} + \varphi_{235}, \mathcal{L}]] + [\varphi_{125}, [\varphi_{134} + \varphi_{234}, \mathcal{L}]].
    \end{align*}
    Expand all commutators using the Jacobi identity for the commutator:
    \begin{align*}
        \widetilde{\RHS} &= -[\mathcal{L}, [\varphi_{123}, \varphi_{145} + \varphi_{245}]] + [\mathcal{L}, [\varphi_{124}, \varphi_{135} + \varphi_{235}]] - [\mathcal{L}, [\varphi_{125}, \varphi_{134} + \varphi_{234}]] \\ &\quad+ [\varphi_{145} + \varphi_{245}, [\varphi_{123}, \mathcal{L}]] -[\varphi_{135} + \varphi_{235}, [\varphi_{124}, \mathcal{L}]] + [\varphi_{134} + \varphi_{234}, [\varphi_{125}, \mathcal{L}]].
    \end{align*}
    Bring all terms of the form $[\mathcal{L}, \ldots]$ to the left-hand side:
    \begin{align*}
        \overline{LHS} &= [\mathcal{L}, [\varphi_{123}, \varphi_{145} + \varphi_{245}] + [\varphi_{125}, \varphi_{134} + \varphi_{234}] - [\varphi_{124}, \varphi_{135} + \varphi_{235}] + [\varphi_{123} + \varphi_{124} + \varphi_{125}, \varphi_{345}]] = \\ &= [\mathcal{L}, [\varphi_{123}, \varphi_{145} + \varphi_{245} + \varphi_{345}] - [\varphi_{124}, \varphi_{135} + \varphi_{235} + \varphi_{435}] + [\varphi_{125}, \varphi_{134} + \varphi_{234} + \varphi_{534}]].
    \end{align*}
    Then the remaining terms on the right are only
    \begin{equation*}
        \overline{RHS} = [\varphi_{145} + \varphi_{245}, [\varphi_{123}, \mathcal{L}]] -[\varphi_{135} + \varphi_{235}, [\varphi_{124}, \mathcal{L}]] + [\varphi_{134} + \varphi_{234}, [\varphi_{125}, \mathcal{L}]].
    \end{equation*}
    Thus the final equation takes the form
    \begin{align*}
        &[\mathcal{L}, [\varphi_{123}, \varphi_{145} + \varphi_{245} + \varphi_{345}] - [\varphi_{124}, \varphi_{135} + \varphi_{235} + \varphi_{435}] + [\varphi_{125}, \varphi_{134} + \varphi_{234} + \varphi_{534}]] \\ &= [\varphi_{145} + \varphi_{245}, [\varphi_{123}, \mathcal{L}]] -[\varphi_{135} + \varphi_{235}, [\varphi_{124}, \mathcal{L}]] + [\varphi_{134} + \varphi_{234}, [\varphi_{125}, \mathcal{L}]].
    \end{align*}
    Note that when writing the conditions for the linear $\varphi$-matrix operation, the same equations appear, but with $\mathcal{L} = L_1 + L_2 + L_3 + L_4 + L_5$. The equation still depends on $\mathcal{L}$ and must hold for all $L$.

\section{From the modified classical Yang--Baxter equation to a 3-Lie algebra}
\label{sec:cyb-fi}

In this Section we generalize our construction to Lie algebras which share common features of $\mathfrak{gl}_N$ and prove that if $r$ satisfies modified classical Yang--Baxter equation then $\varphi(r)$ gives us 3-Lie algebra structure (see Proposition \ref{n-lie-algebra} for connection between Nambu brackets and $n$-Lie algebra structure). 

\subsection{Construction of 3-Lie algebra}
Let $\mathfrak{g}$ be a finite-dimensional real reductive Lie algebra with a nondegenerate symmetric invariant bilinear form $\langle\cdot,\cdot\rangle$, so that
\[
    \langle[x,y],w\rangle=\langle x,[y,w]\rangle.
\]
Choose $z\in Z(\mathfrak{g})$ and set $\tau(x)=\langle z,x\rangle$. The element $z$ plays the role of the identity matrix in the matrix construction. We identify $\mathfrak{g}$ with its dual by the given form and write $\ell_x(L)=\langle x,L\rangle$.

For $r=\sum_\alpha a_\alpha\otimes b_\alpha\in\wedge^2\mathfrak{g}$, use the convention
\begin{align*}
    \operatorname{CYB}(r)
    =\sum_{\alpha,\beta}\bigl(
    [a_\alpha,a_\beta]\otimes b_\alpha\otimes b_\beta
      +a_\alpha\otimes[b_\alpha,a_\beta]\otimes b_\beta+a_\alpha\otimes a_\beta\otimes[b_\alpha,b_\beta]\bigr).
\end{align*}
Also we say that $r \in \wedge^2 \mathfrak{g}$ satisfies modified classical Yang--Baxter equation if $\operatorname{CYB}(r) \in (\wedge^3 \mathfrak{g})^{\mathfrak{g}}$. We would refer to it as $\operatorname{mCYBE}$.
Define the alternating tensor
\begin{equation}\label{eq:central-phi}
    \varphi=r\wedge z
    =\sum_\alpha\bigl(a_\alpha\otimes b_\alpha\otimes z
      -a_\alpha\otimes z\otimes b_\alpha
      +z\otimes a_\alpha\otimes b_\alpha\bigr).
\end{equation}
Let $\ad_L^{(k)}$ denote the sum of the adjoint actions of $L$ on the $k$ tensor factors. The linear ternary operation associated with $\varphi$ is defined by
\begin{equation}\label{eq:quadratic-ternary}
    \{\ell_x,\ell_y,\ell_w\}(L)
    =\bigl\langle-\ad_L^{(3)}\varphi,\,
                   x\otimes y\otimes w\bigr\rangle,
\end{equation}
where the pairing on tensor powers is the product pairing.

\begin{thm}\label{thm:cyb-fi}
    Suppose $r\in\wedge^2\mathfrak{g}$ satisfies $\operatorname{mCYBE}$ and $z\in Z(\mathfrak{g})$. Then \eqref{eq:quadratic-ternary} satisfies the fundamental identity on linear functions. Under the identification by $\langle\cdot,\cdot\rangle$, the resulting 3-Lie bracket on $\mathfrak{g}$ is
    \begin{equation}\label{eq:central-3-lie}
        T(x,y,w)=\langle z,x\rangle B(y,w)
                +\langle z,y\rangle B(w,x)
                +\langle z,w\rangle B(x,y),
    \end{equation}
    where
    \begin{equation}\label{eq:cyb-R}
        R(x)=\sum_\alpha a_\alpha\langle b_\alpha,x\rangle,
        \qquad B(x,y)=[R(x),y]+[x,R(y)].
    \end{equation}
\end{thm}

We use the following result of Arnlind, Kitouni, Makhlouf and Silvestrov, which is the case $n=2$ of their general construction of induced $(n+1)$-Lie algebras.

\begin{lemma}[Theorem 1 \cite{AKMS2013}]\label{lem:trace-3-lie}
    Let $(\mathfrak{a},B)$ be a Lie algebra and let $\tau:\mathfrak{a}\to\mathbb{R}$ be a linear functional such that $\tau(B(x,y))=0$ for all $x,y\in\mathfrak{a}$. Then
    \begin{equation}\label{eq:trace-ternary}
        T(x,y,z)=\tau(x)B(y,z)+\tau(y)B(z,x)+\tau(z)B(x,y)
    \end{equation}
    defines a 3-Lie algebra on $\mathfrak{a}$.
\end{lemma}

\begin{proof}[Proof of Theorem \ref{thm:cyb-fi}]
    According to \cite{STS} from $\operatorname{mCYBE}$ it follows that $B$ is a Lie bracket.
    Centrality of $z$ and invariance imply
    \[
        \tau([u,v])=\langle z,[u,v]\rangle
                   =\langle[z,u],v\rangle=0.
    \]
    Both terms in $B(x,y)$ are Lie brackets in $\mathfrak{g}$, so $\tau(B(x,y))=0$. Lemma \ref{lem:trace-3-lie} therefore proves that \eqref{eq:central-3-lie} is a 3-Lie bracket.

    It remains to identify this bracket with \eqref{eq:quadratic-ternary}. Invariance and skew-symmetry give
    \begin{align*}
        \bigl\langle-\ad_L^{(2)}r,x\otimes y\bigr\rangle
        &=\sum_\alpha\bigl(
          \langle[a_\alpha,L],x\rangle\langle b_\alpha,y\rangle
          +\langle a_\alpha,x\rangle\langle[b_\alpha,L],y\rangle
          \bigr)\\
        &=\langle[R(x),y]+[x,R(y)],L\rangle.
    \end{align*}
    Since $\ad_L z=0$, we also have
    \[
        -\ad_L^{(3)}(r\wedge z)=(-\ad_L^{(2)}r)\wedge z.
    \]
    Pairing with $x\otimes y\otimes w$ now gives
    \begin{equation}\label{eq:cyb-ternary-identification}
        \{\ell_x,\ell_y,\ell_w\}=\ell_{T(x,y,w)},
    \end{equation}
    with $T$ as in \eqref{eq:central-3-lie}. This proves the assertion on linear functions.
\end{proof}

\begin{rem}
    The proof applies to any finite-dimensional Lie algebra with a nondegenerate symmetric invariant form and a chosen central element. In particular, reductivity is sufficient. The form identifies the possible central elements with the linear functionals annihilating $[\mathfrak{g},\mathfrak{g}]$, since
    \[
        Z(\mathfrak{g})=[\mathfrak{g},\mathfrak{g}]^\perp.
    \]
    If $\mathfrak{g}$ is semisimple, this construction gives $z=0$ and hence the zero ternary bracket. We have already seen it in the special case of $\mathfrak{sl}_2$ (Lemma \ref{sl-restriction}).
\end{rem}

\begin{cor}\label{cor:matrix-cyb-fi}
    For $\mathfrak{g}=\mathfrak{gl}_N$, $\langle X,Y\rangle=\Tr(XY)$ and $z=I$, Theorem \ref{thm:cyb-fi} gives a 3-Lie algebra on matrices whenever $r$ satisfies $\operatorname{mCYBE}$. Its corresponding operation on linear functions is
    \[
        \{L_1,L_2,L_3\}
        =[r_{12}-r_{13}+r_{23},L_1+L_2+L_3].
    \]
\end{cor}
\begin{proof}
    Here $\tau(X)=\Tr X$ and
    $R(X)=\Tr_2(r_{12}(I\otimes X))$. The tensor \eqref{eq:central-phi} becomes $r_{12}-r_{13}+r_{23}$, and $-\ad_L^{(3)}\varphi=[\varphi,L_1+L_2+L_3]$. Thus \eqref{eq:quadratic-ternary} is precisely the matrix operation which satisfies linear part of the fundamental identity.
\end{proof}

\subsection{A rank criterion for the fundamental identity}

Now the question is whether trivectors constructed from the solution of $\operatorname{mCYBE}$ are decomposable.

\begin{prop}\label{prop:rank-nambu}
    Under the assumptions of Theorem \ref{thm:cyb-fi}, suppose $z\ne0$ and choose $e\in\mathfrak{g}$ such that $\langle z,e\rangle=1$. Let $\pi$ be the linear trivector defined by \eqref{eq:quadratic-ternary}, and define
    \[
        Q_e(df,dg)=\{f,g,\ell_e\}.
    \]
    Then $Q_e$ is a linear Poisson bivector. The following conditions are equivalent:
    \begin{enumerate}
        \item $\pi$ is a Nambu bracket,
        \item $\pi(L)$ is decomposable at every point $L$ (where it is non-zero),
        \item $\operatorname{rank}Q_e(L)\le2$ at every point $L$.
    \end{enumerate}
\end{prop}
\begin{proof}
    Fixing $e$ in the 3-Lie bracket gives a Lie bracket, so $Q_e$ satisfies the Jacobi identity on linear functions. Since $Q_e$ is linear, it is a Poisson bivector on all smooth functions (Proposition \ref{not-deco}).

    Regard $z$ as a constant vector field on $\mathfrak{g}$ and write
    \[
        \beta_L=-\ad_L^{(2)}r,\qquad \pi(L)=z\wedge\beta_L.
    \]
    Put $\alpha=d\ell_e$, so $\alpha(z)=1$ and contract $\alpha$ with $\pi$. In this Section we use the standard convention of contracting a covector into the first slot of a multivector. For a trivector, contraction in the third slot gives the same bivector by cyclic symmetry.
    \[
        Q_e=\iota_\alpha\pi=\beta-z\wedge\iota_\alpha\beta,
        \qquad \iota_\alpha Q_e=0,\qquad \pi=z\wedge Q_e.
    \]
    Thus $Q_e(L)\in\wedge^2\ker\alpha$, and wedge multiplication by $z$ identifies this space with $z\wedge\wedge^2\ker\alpha$. It follows that $\pi(L)$ is decomposable if and only if $Q_e(L)$ is decomposable, equivalently if and only if its rank is at most two. The same identification gives $\pi(L)=0$ if and only if $Q_e(L)=0$.

    By Theorem \ref{thm:cyb-fi}, $\pi$ already satisfies the fundamental identity on linear functions. The decomposability criterion in Theorem \ref{thm-vaisman} therefore proves the equivalence with the full Nambu property.
\end{proof}

\begin{cor}\label{cor:matrix-rank-nambu}
    Let $r\in\wedge^2\mathfrak{gl}_N$ satisfy $\operatorname{mCYBE}$. The linear ternary operation \eqref{r-matrx-nambu} is a Nambu bracket if and only if the Poisson bracket
    \[
        \{f,g\}_1=\{f,g,\Tr L\}
    \]
    has rank at most two everywhere. If, in addition, $\Tr_2(r_{12})=0$, this is equivalent to the ordinary linear $r$-matrix Poisson bracket
    \[
        \{L_1,L_2\}_{lin}=[r_{12},L_1+L_2]
    \]
    having rank at most two everywhere.
\end{cor}
\begin{proof}
    Apply Proposition \ref{prop:rank-nambu} with $z=I$ and $e=I/N$. Then $\{\cdot,\cdot\}_1=NQ_e$, which has the same rank as $Q_e$. Under the partial-trace condition, direct contraction gives
    \[
        \{L_1,L_2\}_1=N[r_{12},L_1+L_2].
    \]
    This is the identity used in Corollary \ref{cor-bracket}.
\end{proof}

The rank condition refines the necessary conditions discussed in Corollary~\ref{weak-FI}.

\section{Geometric examples and quasi-Hopf twists}\label{sec:principal-examples}

Now we introduce some examples of ternary operations motivated by deformation of (co)associativity (see Section \ref{sec:nonassoc-rtt} for the introduction). We should note that the fundamental identity is not generally guaranteed in these constructions.

\subsection{Heisenberg cochain twist}

Consider the Heisenberg Lie algebra
$\mathfrak h_3=\langle X,Y,Z\rangle$, where
$[X,Y]=Z$ and $Z$ is central, acting on functions on $\mathbb R^3$ by
\begin{equation}\label{heis-fields}
 X=\partial_x,\qquad Y=\partial_y+x\partial_z,
 \qquad Z=\partial_z.
\end{equation}
Take the cochain
\begin{equation}\label{heis-twist}
 F=\exp\left[-\frac{\hbar}{2}
 (X\otimes Y-Y\otimes X)\right].
\end{equation}
This is not a Hopf 2-cocycle, nevertheless it produces a
quasi-Hopf twist of $U(\mathfrak h_3)[[\hbar]]$.
Expanding star product \eqref{twist-star} with the given twist has first-order bivector
\begin{equation}\label{heis-pi}
 \pi=X\wedge Y
 =\partial_x\wedge\partial_y
   +x\partial_x\wedge\partial_z.
\end{equation}
In particular,
\begin{equation}\label{heis-coordinates}
 \{x,y\}=1,\qquad \{x,z\}=x,
 \qquad \{y,z\}=0.
\end{equation}
A direct computation of Jacobiator gives us the canonical Nambu bracket \eqref{canonical-nambu} on $\mathbb{R}^3$
\begin{equation}\label{heis-jacobiator}
 J_\pi(f,g,h)
 =\det\frac{\partial(f,g,h)}{\partial(x,y,z)}.
\end{equation}

\subsection{Jacobiators along three-dimensional leaves}

When we can definitely say that Jacobiator is a Nambu trivector?

\begin{thm}\label{thm:leafwise-jacobiator}
Let $D\subset TM$ be an integrable distribution of rank three on a smooth
manifold $M$, and let $\pi\in\Gamma(\wedge^2 D)$. Then the Jacobiator
$J_\pi=\tfrac12[[\pi,\pi]]$ defines a Nambu bracket of order three
on $C^\infty(M)$.
\end{thm}

\begin{proof}
Involutivity of $D$ implies that the Schouten bracket of multivector fields
tangent to $D$ is again tangent to $D$. Hence $J_\pi$ is a section of
$\wedge^3 D$. Since $D$ has rank three, $J_\pi$ is locally decomposable,
and wherever it is nonzero its characteristic distribution is exactly $D$.
The criterion of \cite[Remark~2.7(i)]{Vaisman99} now gives the fundamental
identity on $M$.
\end{proof}

\begin{ex}
On $\mathbb R^4$ with coordinates $(x,y,z,t)$, take $D=\ker dt$ and
\[
 \pi=\partial_x\wedge\partial_y+tx\,\partial_x\wedge\partial_z.
\]
The only nonzero basic brackets are $\{x,y\}=1$ and $\{x,z\}=tx$. Thus $J_\pi(x,y,z)=\{y,-tx\}=t$, and every
bracket involving $t$ vanishes. The induced ternary bracket is
\[
 J_\pi(f,g,h)=t\,\det\frac{\partial(f,g,h)}{\partial(x,y,z)}.
\]
Its three-dimensional leaves are the hyperplanes $t=c$ for $c\ne0$.
At $t=0$ the Jacobiator vanishes.
\end{ex}

\begin{rem}
On a quasi-Poisson $G$-manifold the Jacobiator is the action trivector
$\rho(\varphi_{\mathrm C})$, where $\varphi_{\mathrm C}$ is the Cartan
trivector
\cite{AlekseevKSM}. If all orbits have dimension at most three, this
Jacobiator is Nambu bracket. Indeed, it is tangent to the orbits and
vanishes on every orbit of dimension less than three. On its nonzero
locus it is decomposable and its characteristic distribution is the
tangent distribution to the three-dimensional orbits, so the same
criterion applies.
\end{rem}

\subsection{$GL_2$ cochain twist}
\label{subsec:gl2-twist-nambu}

In $U_\hbar(\mathfrak{gl}_2)$ take the standard $\mathfrak{sl}_2$
generators $H,E,F$ with $K=e^{\hbar H}$ and coproduct
\[
 \Delta(E)=E\otimes1+K\otimes E,\qquad
 \Delta(F)=F\otimes K^{-1}+1\otimes F.
\]
Let $\star$ be the product obtained by twisting the standard RTT
multiplication on $\mathcal O_\hbar(GL_2)$ via \eqref{twist-star}
with the cochain
\begin{equation}\label{gl2tw-cochain}
 \mathcal F=\exp\!\left[-\tfrac{\hbar}{2}(E\otimes F-F\otimes E)\right].
\end{equation}
The bivector $\beta$ encodes the binary bracket given by the
first-order commutator of this product, as in \eqref{binary-classical}:
\[
 u\star v-v\star u=\hbar\,\beta(du,dv)+O(\hbar^2).
\]
To compute $\beta$, write $X^L(g)=gX$ and $X^R(g)=Xg$.
With our RTT normalization the untwisted first-order commutator has
bivector $\pi_0=E^R\wedge F^R-E^L\wedge F^L$.
We use the action on the coordinate algebra for which $X$ acts
modulo $\hbar$ by the vector field $X^L$. Thus the term linear in
$\hbar$ in $\mathcal F^{-1}$ contributes
$\tfrac12(E^Lu\,F^Lv-F^Lu\,E^Lv)$ to the first-order product
coefficient. Taking the commutator gives
\[
 \beta=\pi_0+E^L\wedge F^L=E^R\wedge F^R.
\]
Put $\varphi=H\wedge E\wedge F$.

\begin{prop}\label{prop:gl2tw-nambu}
The Jacobiator of $\beta$ is $-\varphi^R$. The trivector
$\Pi_N=\varphi^R$ defines a Nambu bracket with the
determinant level sets as its leaves. In matrix coordinates
$g=\left(\begin{smallmatrix}a&b\\c&d\end{smallmatrix}\right)$,
$D=ad-bc$, its values on the coordinate functions are
\begin{equation}\label{gl2tw-coordinate-bracket}
 \begin{aligned}
  \{a,b,c\}_N &= aD, & \{a,b,d\}_N &= bD,\\
  \{a,c,d\}_N &=-cD, & \{b,c,d\}_N &=-dD.
 \end{aligned}
\end{equation}
For arbitrary functions $u,v,w$ this gives
\begin{equation}\label{gl2tw-jacobian}
 \{u,v,w\}_N=-D\,\frac{\partial(D,u,v,w)}{\partial(a,b,c,d)}.
\end{equation}
\end{prop}
\begin{proof}
Since $[E^R,F^R]=-H^R$, the Jacobiator in the convention of
\eqref{jacobi-associator} is $-H^R\wedge E^R\wedge F^R$.
These three fields span $\ker dD$, the tangent distribution to
the $SL_2$-orbits. Thus Theorem~\ref{thm:leafwise-jacobiator}
applies to $\beta$, and its Jacobiator is nonzero everywhere.
For the coordinate formula, let $\Omega=da\wedge db\wedge dc\wedge dd$.
At the identity $\iota_{\varphi^R}\Omega=dD$. Under right
translation by $g$, $\Omega$ scales by $D(g)^2$ and $D$ scales
by $D(g)$, so $\iota_{\varphi^R}\Omega=D\,dD$ throughout $GL_2$.
This gives \eqref{gl2tw-jacobian}.
\end{proof}

Finally, \eqref{jacobi-associator} identifies the leading alternating
associator of the twisted product with the same bracket:
\begin{equation}\label{gl2tw-associator-nambu}
 6\,\operatorname{Alt}\operatorname{As}_2(u,v,w)
 =-J_\beta(u,v,w)=\{u,v,w\}_N.
\end{equation}

\section{Conclusions and open problems}\label{sec:conclusions}

The $\varphi$-matrix construction extends the commutator formulas for
$r$-matrix Poisson brackets to ternary operations on matrix space.
It gives explicit formulas in matrix coordinates and allows us to ask
when the resulting trivectors satisfy the fundamental identity.
For the linear operation associated with
$\varphi=r_{12}-r_{13}+r_{23}$, we obtain a criterion under the
assumption that $r$ satisfies the modified classical Yang--Baxter
equation. The ternary bracket on linear functions then defines a
$3$-Lie algebra on $\mathfrak{gl}_N^*$. Moreover, the full fundamental
identity holds if and only if the Poisson bracket obtained by fixing
the trace,
\[
    \{f,g\}_1=\{f,g,\operatorname{Tr}L\},
\]
has rank at most two at every point.

This rank bound expresses the requirement that the linear trivector be
pointwise decomposable. It also suggests a question about applications:
under what weaker hypotheses on the trivector can one still obtain
compatible Poisson brackets by fixing one argument? We would like to
understand which choices of that argument yield such families and how
they can be used to construct functions in involution.

A second source of ternary operations is the associativity defect of a
deformed product. In the RTT setting we examined the constraints imposed
by the quasi-associative Yang--Baxter equation and, for a universal
coassociator, by the pentagon equation. Whether these constraints imply
the fundamental identity under additional hypotheses remains open.
For the Jacobiator of a bivector, we obtained a geometric sufficient
condition: tangency of the bivector to an integrable distribution of
rank three. Applied to a quasi-Hopf cochain twist on $GL_2$, this
criterion gives a Nambu bracket whose characteristic leaves are the
level sets of the determinant.

\paragraph{Acknowledgements}{The work of D. T. was carried out within the framework of a development program for the Regional Scientific and Educational Center of the P.~G.~Demidov Yaroslavl State University with financial support from the Ministry of Science and Higher Education of the Russian Federation (Agreement on provision of subsidy from the federal budget No. 075-02-2026-1331). The work of M.Ch. and D.T. was partially supported by a grant Leader Mathematics \#25-7-1-4-1 of the Basis foundation.}

\section*{Declaration of generative AI use}
OpenAI's ChatGPT was used for exploratory discussions of
mathematical arguments, assistance with symbolic computations,
and preparation of the \LaTeX{} manuscript.
The authors independently take responsibility for all mathematical
statements, proofs, references and the final text.

\bibliographystyle{unsrt}
\bibliography{main}

\end{document}